\documentclass[10.5pt, twoside, leqno]{article}
\usepackage{xcolor,gensymb}     
\usepackage{amsmath,amsthm}
\usepackage{amssymb}
\usepackage{amscd}   
\usepackage{graphicx}
\usepackage[T1]{fontenc}   

\usepackage{amsmath,amsthm}
\usepackage{amssymb}

\usepackage{enumitem}

\usepackage{graphicx}

\usepackage[T1]{fontenc}

\usepackage{color}
\usepackage[normalem]{ulem}
\newcommand{\ev}{{\rm Eval}}  
\newcommand{\inte}{{\rm Prim}} 
\newcommand{\deri}{{\rm Deri}}   
\newcommand{\id}{{\rm Id}} 
\newcommand{\card}{{\rm Card}}

\newtheorem{theorem}{Theorem}[section]
\newtheorem{lemma}[theorem]{Lemma}
\newtheorem{corollary}[theorem]{Corollary}
\newtheorem{proposition}[theorem]{Proposition}

\theoremstyle{definition}
\newtheorem{definition}[theorem]{Definition}
\newtheorem{remark}[theorem]{Remark}
\newtheorem{example}[theorem]{Example}

\numberwithin{equation}{section}
\newtheorem{facts}[theorem]{Facts}
\newtheorem{notation}[theorem]{Notation}

\newcommand{\C}{\mathbb{C}} 
\newcommand{\R}{\mathbb{R}}
\newcommand{\ru}{{\R}}
\newcommand{\Q}{\mathbb{Q}}

\newcommand{\Z}{\mathbb{Z}} 
\newcommand{\qu}{{\Q}}
\newcommand{\zu}{{\Z}}
\newcommand{\N}{\mathbb{N}}
\newcommand{\enn}{\N}

\newcommand{\pu}{{\mathbb{P}}}

\def\2{I\hspace{-.1em}I}

\newcommand{\consta}{{{\rho} m\log(2)+{\rho} \left(\log({{\rho}} m+1)+{\rho} m\log\left(\dfrac{{\rho} m+1}{{\rho} m}\right)\right)}}
\newcommand{\const}{{{{\rho}} m\log(2)+{{\rho}} \left(\log({{\rho}} m+1)+{{\rho}} m\log\left(\dfrac{{{\rho}} m+1}{{{\rho}} m}\right)\right)}}

\renewcommand{\proofname}{Proof.}
\begin{document}

\title{Linear independence criteria \\for generalized polylogarithms with distinct shifts} 
\author{Sinnou~David, Noriko~Hirata-Kohno and Makoto~Kawashima}
\date{}

\maketitle

\begin{abstract}
For a given rational number $x$ and an integer $s\geq 1$, let us consider a generalized polylogarithmic function, often called the Lerch function, defined by
$$\Phi_{s}(x,z)= \sum_{k=0}^{\infty}\frac{z^{k+1}}{(k+x+1)^s}\enspace.$$
We prove the linear independence over any number field $K$ of the numbers $1$ and $\Phi_{s_j}(x_j,\alpha_i)$ 
with \emph{any choice of distinct shifts} $x_1,\ldots, x_d$ with $0\le x_1<\ldots<x_d<1$, as well as any choice of 
depths $1\leq s_1\leq r_1,\ldots, 1\leq s_d\leq r_d$, 
at distinct algebraic numbers $\alpha_1,\ldots,\alpha_m\in K$ subject to a metric condition.
As is usual in the theory, the points $\alpha_i$ need to be chosen sufficiently close to zero with respect to
a given fixed place $v_0$ of $K$, Archimedean or finite.

This is the first linear independence result with distinct shifts $x_1, \ldots, x_d$ that allows values at different points for generalized polylogarithmic functions. 
Previous criteria were only for the functions with one fixed shift or at one point.

Further, we establish another linear independence criterion for values of the generalized polylogarithmic function \emph{with cyclic coefficients}.
Let $q\geq 1$ be an integer and $\boldsymbol{a}=(a_1,\ldots, a_q)\in K^q$ be a $q$-tuple whose coordinates supposed to be cyclic with the period $q$.
Consider the generalized polylogarithmc function with coefficients
$$\Phi_{\boldsymbol{a},s}(x,z)=
\sum_{k=0}^{\infty}\frac{a_{k+1\bmod(q)}\cdot z^{k+1}}{(k+x+1)^s}\enspace.$$ 
Under suitable condition, we show that the values of these functions are linearly independent over $K$.

Our key tool is a new non-vanishing property for a generalized Wronskian of Hermite type
associated to our explicit constructions of Pad\'e approximants for this family of generalized polylogarithmic function.

\end{abstract}


\section{Introduction}

Pad\'e approximation ({\it confer}~\cite{Pade1,Pade2})
is a classical method in diophantine problems and has been specially efficient in proving  irrationality  of special values of locally holomorphic functions  as well as providing sharp estimates for the associated irrationality measures in question. 
With a rich history and contributions from a large number of authors after the pioneering work of Ch. Hermite and H. Pad\'e, the method has been extensively studied.

For values of polylogarithmic functions or generalized ones, E.~M.~Niki\v{s}in  in \cite{N,N2},
A.~I.~Galochkin in  \cite{G1} and  \cite{G2}, G.~V.~Chudnovsky in \cite{ch11}, K.~V$\ddot{\text{a}}$$\ddot{\text{a}}$n$\ddot{\text{a}}$nen in \cite{Va} and W. Zudi\-lin in \cite{Z} gave several linear independence criteria, either over the  field of rational numbers, 
or quadratic imaginary fields, albeit with only one fixed shift or at one point.
Niki\v{s}in succeeded in proving
the linear independence of polylogarithms in \cite{N} by using Pad\'e approximation of  type I,  and considered a variant of type II  in a more general setup \cite{N2}.
Nikisin's results were with the single shift equal to zero, corresponding $r=s, m=1$ of our case, over the base field being the rational number field. 
Our approach is inspired by these contributions, {{together with that by G.~Rhin and P.~Toffin showing Pad\'e approximation  of type II for logarithms at different points \cite{R-T}.}}

Considering one fixed shift but allowing all other parameters to be free, that was also one main purpose of our previous papers \cite{DHK2,DHK3}.
See these articles for further reference on earlier works and a historical survey in this context. 
Note that our linear independence result of values for the functions in the current paper does \emph{not} automatically follow from the criterion of G.~V.~Chudnovsky in  \cite{ch11}.

Our aim in the present article is to {\it explore} a largely unchartered territory of families of generalized polylogarithmic functions with total freedom on the shifts, on the points, on cyclic coefficients,
and parameters considered. 
Note that M. Hata \cite{Ha} dealt with the special case of evaluation at a single rational point and proved the linear independence over the rational number field of values of generalized polylogarithmic functions with distinct shifts, at one rational point. 
Our construction of Pad\'e approximants was applied in the forthcoming article \cite{KawashimaPoels1} where the
authors established, in {{addition}}, 
an effective version of the Poincar\'e-Perron theorem to prove new irrationality measures.

It is a classical remark that the Riemann $\zeta$ function is expressed as $\zeta(s)=\Phi_{s}(0,1)$. 
Similarly, the Dirichlet $L$-function is viewed as values at $z=1$ of a generalized polylogarithmic function with periodic coefficients and the shift $x=0$, namely for
any given sequence $\boldsymbol a=(a_{k+1})_{k\geq 0}$ of algebraic numbers, we consider {the twisted polylogarithmic function by Hadamard product $\star $ by $\boldsymbol{a}$}:
$${\boldsymbol a}\star \Phi_{s}(x,z)=\sum_{k=0}^{\infty}\frac{a_{k+1} \cdot z^{k+1}}{(k+x+1)^s}\enspace.$$
We may come back to this point at a later stage.

The special case when the sequence ${\boldsymbol a}$ is actually periodic (say of period $q$) is of particular interest (and more so if $\boldsymbol{a}$ is an arithmetic function such as a Dirichlet character). The values at $z=1$ of such functions also encompass the Dirichlet $L$-functions as remarked above.
We refer to previous important works due to T.~Rivoal \cite{Ri}, R.~Marcovecchio \cite{marc}, S.~Fischler \cite{Fis2}, 
Fischler-J.~Sprang-Zudilin \cite{FSZ}, {{Li Lai-Pin Yu \cite{L-Y}}}, where 
lower bounds for the dimension of the vector space are proven, for those spanned by polylogarithms or special values of the Riemann zeta functions, or those of the $L$-functions with Dirichlet characters 
({\it see also} Nishimoto's work \cite{Nishi}). {{Hence neither the irrationality nor the linear independence of the values follows from their works. On the other hand, in our work, the point
$z$ needs to be taken close to the origin, thus for the values of the Riemann zeta function or the Dirichlet $L$-function, our method does not directly give the linear independence result.}}

\section{Notations and Main results}
{Throughout the article, we denote by $\N$
the set of strictly positive integers.
Let $K$ be a number field.} We denote the set of places of $K$  by ${{\mathfrak{M}}}_K$ (respectively by ${\mathfrak{M}}^{\infty}_K$ for archimedean places, by ${{\mathfrak{M}}}^{f}_K$
for finite places).
For $v\in {{\mathfrak{M}}}_K$, we denote by $K_v$ the completion of $K$ with respect to $v$, and by $\mathbb C_v$
the completion of the algebraic closure $\overline{K_v}$ of $K_v$.

For $v\in {{\mathfrak{M}}}_K$, we define the normalized absolute value $| \cdot|_v$ as follows~:
$$
|p|_v=p^{-\tfrac{[K_v:\Q_p]}{[K:\Q]}} \ \text{if} \ v\in{{\mathfrak{M}}}^{f}_K \ \text{and} \ v|p\enspace, \hspace{45pt}
|x|_v=|\iota_v( x)|^{\tfrac{[K_v:\R]}{[K:\Q]}} \ \text{if} \ v\in {{\mathfrak{M}}}^{\infty}_K\enspace,
$$
where $p$ is a rational prime and $\iota_v$ the embedding $K\hookrightarrow \C$ corresponding to $v$. 
The norm $\Vert\cdot\Vert_v$ denotes the norm of the supremum in $K_v^n$.
{With these normalizations, the product formula reads}
\begin{align*} 
\prod_{v\in {{\mathfrak{M}}}_K} |\xi|_v=1 \ \text{for} \ \xi \in K\setminus\{0\}\enspace.
\end{align*}

Let $m\in \N$ and ${\boldsymbol{\beta}:=(\beta_0,{\ldots},\beta_m) \in\pu_m( K)}$.  
We define the absolute Weil height of $\boldsymbol{\beta}$ by
\begin{align*}
&{{\mathrm{H}}(\boldsymbol{\beta})=\prod_{v\in {{\mathfrak{M}}}_K} \max\{ |\beta_0|_v,\ldots,|\beta_m|_v\}\enspace,}
\end{align*}
and denote the logarithmic absolute Weil height by ${\rm{h}}(\boldsymbol{\beta})={\rm{log}}\, \mathrm{H}(\boldsymbol{\beta})$. 
For each $v\in \mathfrak{M}_K$, put  ${\rm{h}}_v(\boldsymbol{\beta})=\log\Vert \boldsymbol{\beta}\Vert_v$. 
Then we have ${\mathrm{h}}(\boldsymbol{\beta})={\displaystyle{\sum_{v\in \mathfrak{M}_K}}}{\mathrm{h}}_v(\boldsymbol{\beta})$.  Finally 
viewing $\beta\in K\hookrightarrow \pu_1(K)$ via the natural embedding $\beta\longmapsto (1:\beta)$, we have
\begin{align*}
&{{\mathrm{H}}({\beta})=\prod_{v\in {{\mathfrak{M}}}_K} \max\{ 1,|\beta|_v\}\enspace.}
\end{align*}
Let $S\subset \overline{\Q}$ be a finite set. Define the {\it denominator} of this finite set $S$ by 
$${\rm{den}}(S)=\min \{1\leq n \in \Z\mid n\alpha~\text{is an algebraic integer for each}~ \alpha \in S \}\enspace.$$

For a place $v\in \mathfrak{M}_K$ and $f(z)=\sum_{k=0}^{\infty}f_k/z^{k+1}\in (1/z)\cdot K[[1/z]]$, 
we denote the embedding corres\-ponding to $v$ by $\iota_v:K \hookrightarrow K_v{{\subset \mathbb{C}_v}}$ and
write $f_v(z)={{\sum_{k=0}^{\infty}}}\iota_v(f_k)/z^{k+1}\in (1/z) \cdot\mathbb{C}_v[[1/z]]$.
We then consider $f_v(z)$ as a $v$-adic analytic function inside its {{disk}} of convergence.

{Consider the generalized polylogarithmic function
$$\Phi_{s}(x,z)= \sum_{k=0}^{\infty}\frac{z^{k+1}}{(k+x+1)^s}\enspace,$$
well defined for rational number $x$ which is not negative integer and for
any integer $s\geq 1$ inside its {{disk}} of convergence with the above conventions. 

We are now in the position to introduce all the notations necessary to state our main result.
Let $m, d, r_1,\ldots, r_d \in \N$ and $\boldsymbol{x}:=(x_1, \ldots, x_d)\in\Q^d$ {{satisfying $x_i-x_j\notin \Z$ for $1\leq i<j\leq d$ where {{we}} may assume $0\leq x_d<\cdots <x_1<1$ without loss of generality.
Let $\boldsymbol{\alpha}:=(\alpha_1, \ldots, \alpha_m)\in (K{{\setminus \{0\}}})^m$ whose coordinates being pairwise 
distinct algebraic numbers. 
Put ${{\rho}}={{\sum_{j=1}^d}}r_j$.  
Let $v_0$ be a place of $K$, $\beta\in K$ with $\Vert \boldsymbol{\alpha} \Vert_{v_0}<|\beta|_{v_0}$. Put $b=\displaystyle{\max_{1\le j \le d}}{\rm{den}}(x_j)$.
For a rational number $x$, let us define $$\mu(x)={\rm{den}}(x){\displaystyle{\prod_{\substack{q:\text{prime} \\ q|{\rm{den}}(x)}}}}q^{1/(q-1)}\enspace,$$ and write
\begin{align*} 
V(\boldsymbol{\alpha},\boldsymbol{x},\beta)&=V_{v_0}(\boldsymbol{\alpha},\boldsymbol{x},\beta)={\rm{log}}\,|\beta|_{v_0}-{{\rho}} m{\mathrm{h}}(\boldsymbol{\alpha},\beta)-{{\rho}} m\log\,\|\boldsymbol{\alpha}\|_{v_0}+{{\rho}} m\log\,\|(\boldsymbol{\alpha},\beta)\|_{v_0}
-\sum_{j=1}^dr_j\log\,\mu(x_j)\\
&-\left[\consta\right]- \max_{j}\, (r_j)\cdot b{{\rho}} m\enspace.
\end{align*}

Under the notations above, we give the following linear independence criterion first in a simple manner, 
for values of the generalized polylogarithmic functions with distinct shifts and several depths at different algebraic points.
\begin{theorem} \label{Lerch}
Assume $V(\boldsymbol{\alpha},\boldsymbol{x},\beta)>0$. Then the ${{\rho}} m+1$ numbers~$:$
$$1,\Phi_{s_j}(x_j,\alpha_i/\beta) \ \ \ (1\le i \le m, \ 1\le j \le d, \ 1\le s_j \le r_j)\enspace,$$
are linearly independent over $K$.
\end{theorem} 

\begin{remark}\label{effectif}  
{Our approach is entirely effective, indeed, it is possible to give effective irrationality measures and
effective linear independence measures. {For details, see Theorem $\ref{Lerch 2}$, the more complete, albeit more technical statement.}
We should also mention that our main theorem does \emph{not} immediately follow from the statement nor the proof in \cite{DHK4}, since the case with distinct shifts
in our current setup cannot be covered by the argument in  \cite{DHK4}.}
\end{remark}

{\begin{remark} \label{rem V>0}
Fix $\boldsymbol{\alpha},\boldsymbol{x}$ and assume
$|\beta|_{v_0}>\max(1,\|\boldsymbol{\alpha}\|_{v_0})$. 
Regarding $\beta$ as a variable, the number 
$${\rm{log}}\,|\beta|_{v_0}-{{\rho}} m{\mathrm{h}}(\boldsymbol{\alpha},\beta)+{{\rho}} m\log\,\|(\boldsymbol{\alpha},\beta)\|_{v_0}=\log\,|\beta|_{v_0}-{{\rho}} m\displaystyle{\sum_{v\neq v_0}}{\rm{h}}_{v}(\boldsymbol{\alpha},\beta)\enspace$$ 
is only the part of $V(\boldsymbol{\alpha},\boldsymbol{x},\beta)$ depending on $\beta$, and other terms do not depend on $\beta$. 
Whenever the number 
$$\log\,|\beta|_{v_0}-{{\rho}} m\sum_{v\neq v_0}{\rm{h}}_{v}(\boldsymbol{\alpha},\beta)\enspace,$$ is sufficiently large, then the condition $V(\boldsymbol{\alpha},\boldsymbol{x},\beta)>0$ is verified. On the analytic side, we have tried to keep computations no too technical, at the cost here and there of optimality.

We should note that the condition $V(\boldsymbol{\alpha},\boldsymbol{x},\beta)>0$ indeed depends on the field $K$ (see the proof of Corollary below).
\end{remark}}

\begin{corollary} \label{SAT}
Let $K$ be an algebraic number field and $v_0$ be a place of $K$.
Under  the same notations as in Theorem $\ref{Lerch}$, there exist infinitely many $\beta\in K$ where the ${{\rho}} m+1$ numbers~$:$ 
$$1,\Phi_{s_j}(x_j,\alpha_i/\beta) \ \ \ (1\le i \le m, \ 1\le j \le d, \ 1\le s_j \le r_j)\enspace,$$
are linearly independent over $K$.
\end{corollary}

\proofname

We fix a finite set of places $\mathfrak{M}$ of $K$ of cardinality, say $h\geq 2$, with $v_0\notin \mathfrak{M}$.
We denote $\mathfrak{M}=\{v_1,\ldots,v_h\}$. 
Let $\varepsilon$ be a positive real number with $\varepsilon<1/2$.
By strong approximation theorem ({\it confer}  \cite[Theorem $1.2$]{Lang}), there exists an element $\beta=\beta(\varepsilon) \in K$ which depends on $\varepsilon$ with 
\begin{align} \label{sat}
|\beta-1|_{v_1}<\varepsilon, \ \ \ |\beta|_{v_i}<\varepsilon \ \ \text{for} \ 2 \le i \le h, \ \ \ |\beta|_v\le 1 \ \ \text{for} \ v\notin \mathfrak{M}\cup\{v_0\} \enspace.
\end{align} 
By the condition $|\beta-1|_{v_1}<\varepsilon$, we have $\beta\neq 0$ and $|\beta|_{v_1}\le 2$.
We use the product formula for $\beta$, to obtain 
\begin{align*}
1=|\beta|_{v_1} \cdot \prod_{i=2}^h|\beta|_{v_i} \cdot |\beta|_{v_0} \cdot \prod_{v \notin \mathfrak{M}\cup \{v_0\}}|\beta|_{v}
  \le 2\varepsilon^{h-1}|\beta|_{v_0}\enspace.
\end{align*}
Thus we conclude 
\begin{align*}
1<\dfrac{1}{2\varepsilon^{h-1}}\le |\beta|_{v_0}\enspace. 
\end{align*}
Using the inequalities above, we have 
\begin{align} \label{ineq epsilon}
\log\,|\beta|_{v_0}-{{\rho}} m\sum_{v\neq v_0}{\rm{h}}_{v}(\boldsymbol{\alpha},\beta)& \ge {\rm{log}}\,|\beta|_{v_0}-{{\rho}} m({\rm{h}}(\boldsymbol{\alpha})+\log(2)) \nonumber \\
                                                                &\ge {\rm{log}}\left(\dfrac{1}{2\varepsilon^{h-1}}\right)-{{\rho}} m({\rm{h}}(\boldsymbol{\alpha})+\log(2))\nonumber\\
&\ge -(h-1){\rm{log}}(\varepsilon)-({{\rho}} m+1)({\rm{h}}(\boldsymbol{\alpha})+\log(2))\enspace.
\end{align}
Take a real number $\varepsilon_1$ with $0<\varepsilon_1<1/2$ to have
\begin{align*}
(h-1){\rm{log}}(\varepsilon_1)&<-{{\rho}} m\log\,\|\boldsymbol{\alpha}\|_{v_0}-\sum_{j=1}^dr_j\log\,\mu(x_j)- \max_{j}\, (r_j)b{{\rho}} m-({{\rho}} m+1)({\rm{h}}(\boldsymbol{\alpha})+\log(2))\\
&-\left[\consta\right] \enspace,
\end{align*}
and $\beta_1:=\beta(\varepsilon_1)$ satisfying the inequalities $(\ref{sat})$ for $\varepsilon_1$.
Then, by the definition of $V(\boldsymbol{\alpha},\boldsymbol{x},\beta_1)$ and $(\ref{ineq epsilon})$, we have $V(\boldsymbol{\alpha},\boldsymbol{x},\beta_1)>0$.
Inductively, we take the sequences $(\varepsilon_k)_k\in \R^{\enn }$ and $(\beta_k)_{k}\in K^{\enn}$ with $\varepsilon_{k+1}=\min(|\beta_k-1|_{v_1},|\beta_k|_{v_i})_{2\le i \le h}$ and $\beta_{k+1}=\beta(\varepsilon_{k+1})$.
Then we have $\cdots<\varepsilon_{k+1}<\varepsilon_k<\cdots<\varepsilon_1$, $\beta_{k_1}\neq \beta_{k_2}$ for any $k_1\neq k_2$ and $V(\boldsymbol{\alpha},\boldsymbol{x},\beta_k)>0$ for any $k$. 
By Theorem $\ref{Lerch}$, for any $\beta_k$, the ${{\rho}} m+1$ numbers~:
$$1,\Phi_{s_j}(x_j,\alpha_i/{{\beta_k}}) \ \ \ (1\le i \le m,\ 1\le j \le d, \ 1\le s_j \le r_j)\enspace,$$
are linearly independent over $K$. \qed

Now we {turn to the situation of \emph{generalized polylogarithmic functions}. It is indeed useful to put them inside the more general class of such functions
of which they only happen to be a special case, and treat altogether the general case.
Let $q\geq 1$ be an integer, ${{b}}(z)\in K[z]$ be a polynomial  of degree $q$  and let finally $w=w(z)$ be a polynomial of degree $\leq q-1$ in $K[z]$. 
Let $f_{{{b}},w}(z)=\sum_{k=0}^{\infty}{b_{w,k}}/{z^{k+1}}$ be the power series expansion of the rational function $w(z)/{{b}}(z)$.}

Then, for $x\in\qu$ {{which is not negative integer}} and for $s\geq 0$, the associated generalized polylogarithmic function with coefficients is defined by~:
$$f_{{{b}},w,{{x}},s}(z)=\sum_{k=0}^{\infty}\frac{b_{w,k}\cdot {z^{-k-1}}}{(k+x+1)^s}\enspace.$$
Note that for $s=0$, $f_{{{b}},w,{{x}},0}$ is just the rational function $f_{{{b}},w}$.
{\begin{example} 
In the special case where ${{b}}(z):=z^q-\alpha^q$ {{and $w(z):=\sum_{l=1}^qa_{l}\alpha^{l}z^{q-l}$}}, the function above has the form below with periodic coefficients~:
$$f_{{{b}},w,{{x}},s}(z)=\sum_{k=0}^{\infty}\frac{a_{k+1}}{(k+x+1)^s}\frac{\alpha^{k+1}}{z^{k+1}} \ \ \ \text{with} \ a_{k+1}=a_{l+1} \ \text{if} \  k\equiv l \bmod(q)\enspace,$$ 
and {we indeed recover the statement of Theorem 2.1
for all such functions as $f_{b,w,x,s} (z)$  with periodic coefficients, 
instead of $ \Phi_{s} (x, 1/z)$.}

\end{example}}
Note that by the definition of $f_{{b},w,x,s}(z)$, we have
\begin{align*}
{{b}}(z)\left(-z\dfrac{d}{dz}+x\right)^{s} {{( f_{{{b}},w,x,s}(z))}}=w(z)\in K[z]\enspace.
\end{align*}
 We now state our result for generalized polylogarithmic functions with periodic coefficients.
{{\begin{theorem} \label{periodiclerch}
Let $m,d\geq 1$. 
Assume that we are given $r_1,\ldots r_d$ integers $\geq 1$, put ${{\rho}}=\sum_{j=1}^dr_j$. Let ${{b}}(z)\in K[z]$ be a polynomial of degree $m$ with simple roots $\alpha_1,\ldots,\alpha_m\in K$.
Let $x_1,\ldots,x_d\in\qu$ such that for $i\neq j$, $x_i-x_j\not\in\zu$. Let moreover $w_0(z),\ldots,w_{m-1}(z)$ be {{$K$-}}linearly independent polynomials in $K[z]$ of degree $\leq m-1$.  
Let finally $v_0$ be a place of $K$ and $\beta\in K$ such that $\Vert \boldsymbol{\alpha}\Vert_{v_0}<\vert\beta \vert_{v_0}$.
 
Assume $V({\boldsymbol\alpha},\boldsymbol{x}, \beta)>0$ where $V({\boldsymbol\alpha},\boldsymbol{x}, \beta)$ is defined in Theorem~$\ref{Lerch}$.
 Then, the ${{\rho}} m+1$ numbers, $1, f_{{{b}},w_i,x_j,s_j}(\beta)$ with $0\leq  i\leq m-1$, $1\leq j\leq d$ and $1\leq s_j\leq r_j$ are linearly independent over $K$.
 \end{theorem}
}}

We now state more precisely a criterion as a corollary for the class of functions of  arithmetic interest\,: the generalized polylogarithmic functions with periodic coefficients.
\begin{corollary} \label{coro 1}
Let $K$ be a number field and $v_0$ be a place of $K$. 
Let $m,d,q_1,\ldots,q_m,r_1,\ldots,r_d\in \N$. Put $
{{\rho}}=\sum_{j=1}^dr_j$.
Let $\boldsymbol{\alpha}:=(\alpha_1,\ldots,\alpha_m)\in (K{{\setminus\{0\}}})^m$ and $\beta\in K{{\setminus\{0\}}}$ such that for $i\neq j$, $\alpha_i\neq \alpha_j$ and $\|\boldsymbol{\alpha}\|_{v_0}>|\beta|_{v_0}$.
Let $x_1,\ldots,x_d\in\qu$ such that for $i\neq j$, $x_i-x_j\not\in\zu$. 
Let moreover $b_i(z):=(z^{q_i}-\alpha^{q_i})$ and $w_{i,0}(z),\ldots,w_{i,q_i-1}(z)$ be {{$K$-}}linearly independent polynomials in $K[z]$ of degree $\leq q_i-1$ for $1\le i \le m$.  
Assume $V({\boldsymbol\alpha}',\boldsymbol{x}, \beta)>0$ where $V({\boldsymbol\alpha}',\boldsymbol{x}, \beta)$ is defined in Theorem~$\ref{Lerch}$ for
$$\boldsymbol{\alpha}':=(\alpha_1,\alpha_1\zeta_{q_1},\ldots,\alpha_1\zeta^{q_1-1}_{q_1}, \ldots, \alpha_m,{{\alpha_m}}\zeta_{q_m},\ldots,\alpha_m\zeta^{q_m-1}_{q_m})\enspace,$$
where $\zeta_{q_i}$ is a primitive $q_i$-th {{root}} of unity. 
Then, the $(q_1+\cdots+q_m){{\rho}}+1$ numbers, $1, f_{b_i,w_{i,l_i},x_j,s_j}(\beta)$ with $1\leq  i\leq m$, $0\le l_i \le q_i-1$, $1\leq j\leq d$ and $1\leq s_j\leq r_j$ are linearly independent over $K$.
\end{corollary}

{{\begin{remark}  
{By applying Theorem~\ref{periodiclerch} with the field
$\Q(\alpha_1,\ldots,\alpha_m,\beta,\zeta_{q_1},\ldots,\zeta_{q_m})$ instead of $\Q(\alpha_1,\ldots,\alpha_m,\beta)\subseteq K$,
this corollary holds.}
This follows thanks to the fact that the crucial non-vanishing property
of the determinant in the proof
does \emph{not} depend on the choice of $w_i(z)$. 
Similar observations can be made (specially in the $p$-adic case) and the linear independence can be achieved over some larger but specific fields, however, 
it should be stressed that these methods do not lead to the linear independence over the algebraic closure, neither the transcendence. \end{remark}}}

{{The present article is organized as follows.  In Section~\ref{preli}, we recall the construction of Pad\'e approximants for the generalized polylogarithmic functions as well as various elementary facts that are going to come in the course of the proof.  Section~\ref{deter} contains the crucial non-vanishing result. Although the basic structure of the proof is very classical and dates back to A.~T.~Vandermonde  (separation of variables), its proof is quite involved and extremely demanding. It is thus divided in various subsections and the main steps are outlined at the beginning of the section (see Subsection~\ref{debutdeter}). Section~\ref{estim} is devoted to the proof of our main results after recalling the analytic estimates pertaining to Pad\'e approximation for the generalized polylogarithmic functions. Once non-vanishing proof is achieved, this part proceeds essentially along the lines of our previous paper devoted to a single shift; see \cite{DHK3}. Indeed, the presence of shifts is not 
significant in analytic estimates: they just marginally impact the value of the function $V(\boldsymbol{\alpha},\boldsymbol{x},\beta)$.
Finally, Section~\ref{ex} is devoted to providing concrete examples of Theorem $\ref{Lerch}$.
}}
\section{Preliminaries}
\label{preli}

\subsection{Pad\'{e} approximants}\label{pade}
In this subsection, we explicitly construct Pad\'{e} approximants of the generalized polylogarithmic functions {{with distinct shifts}}. {These can be found in our previous paper~\cite{DHK3}, and are reproduced here for the convenience of the reader}.
First we recall the definition of Pad\'{e} approximants of formal Laurent series.  We denote by $K$ 
a field of characteristic $0$.  
We define the order function at $z=\infty$, ${\rm{ord}}_{\infty}$, by
\begin{align*}
{\rm{ord}}_{\infty}:K((1/z))\rightarrow \Z\cup \{\infty\}; \ \ \sum_{k}{\dfrac{a_k}{z^k}}\mapsto \min\{k\in \Z\mid a_k\neq 0\}\enspace.
\end{align*}  
We first recall without proof the following elementary fact~:
\begin{lemma} \label{padepreliminaries}
Let $r$ be a positive integer, $f_1(z),\ldots,f_r(z)\in (1/z)\cdot K[[1/z]]$ and $\boldsymbol{n}:=(n_1,\ldots,n_r)\in \N^{r}$.
Put $N=\sum_{i=1}^rn_i$.
Let $M$ be a positive integer  with $M\ge N$. Then{,} there exists a family of polynomials 
$(P_0(z),P_{1}(z),\ldots,P_r(z))\in K[z]^{r+1}\setminus\{\bold{0}\}$ satisfying the following conditions~$:$
\begin{align*} 
&({i}) \ {\rm{deg}}\, P_{0}(z)\le M\enspace,\\
&({ii}) \ {\rm{ord}}_{\infty} \, (P_{0}(z)f_j(z)-P_j(z))\ge n_j+1 \ \text{for} \ 1\le j \le r\enspace.
\end{align*}
\end{lemma}
\begin{definition}
We say that a $(r+1)$-tuple of polynomials $(P_0(z),P_{1}(z),\ldots,P_r(z)) \in K[z]^{r+1}$ satisfying the properties $({i})$ and $({{ii}})$ is  a weight $\boldsymbol{n}$ and degree $M$ Pad\'{e} type {II system of} approximants of $(f_1,\ldots,f_r)$.
For such approximants $(P_0(z),P_{1}(z),\ldots,P_r(z))$, of $(f_1,\ldots,f_r)$, 
we call the {$(r+1)$-tuple} of formal Laurent series $(P_{0}(z)f_j(z)-P_{j}(z))_{1\le j \le r}$, {\it id est} the remainders, as weight $\boldsymbol{n}$ degree $M$ Pad\'{e} type approximations of $(f_1,\ldots,f_r)$.
\end{definition}
\bigskip

In the following, we fix $x\in K$ {{which is not negative integer}}.
We now introduce notations for formal primitive, derivation, and evaluation maps.
\begin{notation} \label{notationderiprim}
\begin{itemize}
\item[$(i)$]  For $\alpha\in K$, we denote by ${\ev}_{\alpha}$ the linear evaluation map $K[t]\longrightarrow K$, $P\longmapsto P(\alpha)$. At a later stage, when several variables are in play and there is a perceived ambiguity on which variable is being specialized, we shall denote the map by $\ev_{t\rightarrow \alpha}$.
\item[$(ii)$] For $P\in K[t]$, we denote by $[P]$ the multiplication by $P$ (the map $Q\longmapsto PQ$).
\item[$(iii)$] We also denote by $\inte_x$ the linear operator $K[t]\longrightarrow K[t]$, defined by $P\longmapsto \frac{1}{t^{1+x}}\kern-1,7pt\int_{0}^{t}{\xi}^xP(\xi)d\xi$ (formal primitive). 
\item[$(iv)$] We denote by $\deri_x$ the derivative map $P\longmapsto t^{-x}\tfrac{d}{dt}(t^{x+1}P(t))$, and for $n\geq 1$, $S_{0,x}={\rm{Id}}$, by $S_{n,x}$ the map taking $t^k$ to $\tfrac{(k+x+1)_n}{n!}t^k$ where 
$(k+x+1)_n:=(k+x+1)\cdots(k+x+n)$ 
(the divided derivative $$P\longmapsto \tfrac{1}{n!}t^{-x}\frac{d^n}{dt^n}(t^{n+x}P)=\tfrac{1}{n!}\left(\tfrac{d}{dt}+x/t\right)^n{(t^n P)},$$ so that $\deri_x=S_{1,x}$.
\item[$(v)$] If $\varphi$ is a $K$-automorphism of a $K$-module $M$ and $k$ an integer, we denote 
$$\varphi^{(k)}=\begin{cases}
\overbrace{\varphi\circ\cdots\circ\varphi}^{k-\text{times}} & \ \text{if} \ k>0\\
{\rm{id}}_M &  \ \text{if} \ k=0\\
\overbrace{\varphi^{-1}\circ\cdots\circ\varphi^{-1}}^{-k-\text{times}} & \ \text{if} \ k<0\enspace.
\end{cases} 
$$
\end{itemize}
\end{notation}
\begin{facts}\label{faitselem}
\begin{itemize} 
\item[$(i)$] The map $\inte_x$ is an isomorphism and its inverse is $\deri_x$ for $x\in K$ which is not negative integer. 
\item[$(ii)$] For any integers $n_1,n_2\geq 0$ and $x_1,x_2\in K\setminus\Z$ {{which are not negative integers}}, 
the divided derivatives $S_{n_1,x_1}$ and $S_{n_2,x_2}$ commute, $S_{n_1,x_1}\circ S_{n_2,x_2}=S_{n_2,x_2}\circ S_{n_1,x_1}$.
\item[$(iii)$] By continuity {with respect to  the natural valuation}, all the above mentioned maps {naturally} extend to $K[[t]]$ .
\end{itemize}
\end{facts}
We now quote properties of $S_{n,x}$. 
\begin{lemma} \label{key1} {\rm{\cite[Lemma $3.6$]{DHK3}}}

Let $n$ be a positive integer, $k$ a non-negative integer and $x\in K$ which is not negative integer, one has the following relations valid in ${\rm{End}}_K(K[[t]])$~$:$
\begin{itemize}
\item[$({{i}})$] $S_{n,x}=\dfrac{1}{n!}S_{1,x}\circ(S_{1,x}+\id)\circ\cdots\circ(S_{1,x}+(n-1)\id)\enspace \mbox{and} \enspace [t^k]\circ S_{1,x}=(S_{1,x}-k\id)\circ[t^k]\enspace.$
\item[$({{ii}})$] There exist rational numbers $\{b_{n,m,l}\} \subset \Q$ 
with $b_{n,m,0}=\tfrac{(-m)_n}{n!}$ and, for every {pair of integers} $0\leq m \leq n$,
$$
[t^m]\circ S_{n,x}=\sum_{l=0}^mb_{n,m,l}S^{(l)}_{1,x}\circ[t^m]\enspace.
$$
\end{itemize} 
\end{lemma}

{
\begin{definition} The Pad\'e kernel map associated to the generalized polylogarithmic function $\Phi_s(x,\alpha/z)$ is defined as 
$$\varphi_{\alpha,x,s}=[\alpha]\circ \ev_{\alpha}\circ\inte_{x}^{(s)}\enspace.$$
\end{definition}}
Remark {that} we have
\begin{align} \label{shift}
\varphi_{\alpha, x,s}\circ S_{1,x}=\varphi_{\alpha,x,s-1}\enspace.
\end{align} 
\begin{lemma} \label{kernel} $(${\it confer}  {\rm{\cite{DHK3}}}, Fact $1$ $(v))$ The kernel of the map $\varphi_{\alpha,x,0}$ is the ideal $(t-\alpha)$.
\end{lemma}

We {first} concentrate on a few elementary properties of the maps $\varphi$ which we regroup here and will be useful for the rest~:
\begin{lemma}\label{prelimiphi}
\begin{itemize}
\item[$({{i}})$] The morphisms {$\varphi_{\alpha,x_j,s_j}\circ\inte_{x_j}^{(l)},\varphi_{\beta,x_{j'},s_{j'}}\circ\inte_{x_{j'}}^{(l')}$} pairwise commute for $1\leq j,j'\leq d$, $1\le s_j \le r_j$, $1\le s_{j'} \le r_{j'}$ and $l,l'\in \Z$.
\item[$({{ii}})$] The operator $\frac{\partial}{\partial\alpha}$ commutes with any {$\varphi_{\beta,x_j,s_j}$}.
\item[$({{iii}})$] For any 
$a(t)\in K[t]$, 
we have {{$\varphi_{\alpha,x_j,s_j}(a(t))=[\alpha]\circ\varphi_{1,x_j,s_j}(a(\alpha t))$}.}
\end{itemize}
\end{lemma}

\proofname

{{The assertion}} $({{i}})$ follows from the definition since both multiplication by a scalar, specialization of one variable or integration with respect to a given variable all pairwise commute.
{{The assertion}} $(ii)$ follows from commutation of integrals with respect to a parameter with differentiation with respect to that parameter.

We turn to $(iii)$. Consider $i=1$. 
Take $k\geq 0$, then 
{\begin{align*}
\varphi_{\alpha,x,s}(t^k)&=[\alpha]\circ\ev_{\alpha}\circ\inte^{(s)}_x(t^k)=[\alpha]\circ\ev_\alpha \left( \frac{t^{k}}{(k+x+1)^s}\right)=\frac{\alpha^{k+1}}{(k+x+1)^s}\enspace.
\end{align*} }
Similarly, {$[\alpha ]\circ\varphi_{1,x,s}((\alpha t)^k)=\alpha^{k+1}/(k+x+1)^s$}. Since both coincide, for the basis $\{t^k\}_{k{{\ge 0}}}$, they coincide on ${{K}}[t]$ by linearity.
This completes the proof of Lemma $\ref{prelimiphi}$. \qed

We {now construct} type I\hspace{-.1em}I Pad\'{e} approximants of {the} generalized polylogarithmic functions {$(\Phi_{s_j}(x_j,\alpha_i/z))_{\substack{1\le i \le m,  \\ 1\le j \le d, 1\le s_j \le r_j}}$.}
Let $l$ be a non-negative integer with $0\le l \le {{\rho}} m$ {(recall  ${{\rho}}=\sum_{j=1}^d r_j$)}. For a {{non-negative integer}} $n$, we define a {tuple} of polynomials. 
{Set
\begin{align}
\label{beta}
&A_{l}({t})={t}^l\prod_{i=1}^m({t}-\alpha_i)^{{{\rho}} n}\enspace,
\\
&P_{l}(z)={\rm{Eval}}_z{\Big(\circ_{j=1}^dS^{(r_j)}_{n,x_j}\left(A_l({t})\right)\Big)\enspace, \label{Qnl}}\\
&{ P_{l,i,j,s_j}(z)=}{{P_{l,i,s_j}(z)=\varphi_{\alpha_i,x_j,s_j}}\left(\dfrac{P_{l}(z)-P_{l}(t)}{z-t}\right)\enspace, \label{Qnlijsj}}
\end{align}}
for $1\le i \le m,  1\le j \le d$. {To ease notations, we have deleted the index $j$ in the definition of $P_{l,i,s_j}$. However, some ambiguity might arise when ${ \alpha_i},x_j,s_j$ are specialized to some fixed values.  To avoid it, we adopt the convention that the integer $s_j$  always belongs to the tagged set of integers  $\{_j\kern-1pt1,_j\kern-1pt2,\ldots,_j\kern-1pt r_{j}\}$}{$\simeq\{(j,1),\ldots,(j,r_j)\}$}. }

Under the above notations, we obtain the following theorem. {Note that these polynomials in fact also depend on $n$, we suppressed the extra {{subscript}} to ease reading}.
\begin{theorem} \label{Pade appro Lerch} 
The tuple of polynomials {$(P_{l}(z),P_{l,i,s_j}(z))_{\substack{1\le i \le m,  \\ 1\le j \le d, 1\le s_j \le r_j}}$ form a weight $(n,\ldots,n)\in \N^{{{\rho}} m}$ and degree ${{\rho}} mn+l$
Pad\'{e} type {II system of} approximants of the generalized polylogarithmic functions $$(\Phi_{s_j}(x_j,\alpha_i/z))_{\substack{1\le i \le m, \\ 1\le j \le d, 1\le s_j \le r_j}}\enspace.$$}
\end{theorem} 

\proofname

{This is a variation of Theorem $2$ in \cite{DHK3}}.
By the definition of $P_{l}(z)$, we have
$$ 
{\rm{deg}}\, P_{l}(z)={{\rho}} mn+l\enspace,
$$ 
so the degree condition is satisfied. We need only to check the valuation condition.

Put {{$R_{l,i,s_j}(z)=P_{l}(z)\Phi_{s_j}(x_j,\alpha_i/z)-P_{l,i,s_j}(z)$.}}
Then 
we obtain 
\begin{align*}
R_{l,i,s_j}(z)=\displaystyle P_{l}(z)\varphi_{\alpha_i,x_j,s_j}\left(\dfrac{1}{z-t}\right)-P_{l,i,s_j}(z) 
                                                       =\displaystyle\varphi_{\alpha_i,x_j,s_j}\left(\dfrac{P_{l}(t)}{z-t}\right)
                                                       =\sum_{k=0}^{\infty}\dfrac{\varphi_{\alpha_i,x_j,s_j}(t^kP_{l}(t))}{z^{k+1}}\enspace. 
                  \end{align*}
                 We only  need to show
                  \begin{align*}
\varphi_{\alpha_i,x_j,s_j}(t^kP_{l}(t))=0 \ \text{for} \ 1\le i \le m, 1\le j \le d, 1\le s_j \le r_j \ \text{and} \ 0\le k \le n-1\enspace.
\end{align*}
By Lemma~\ref{key1},  $[t^k]\circ S_{n,x_j}= U(S_{1,x_j})\circ [t^k]$ where $U\in \mathbb{Q}[X]$ is a degree $n$ and valuation $\geq 1$ polynomial, provided that $k\leq n-1$. By iteration, $[t^k]\circ S_{n,x_j}^{(r_j)}=S_{1,x_j}^{(r_j)}\circ V(S_{1,x_j})\circ[t^k]$ where $V\in \mathbb{Q}[X]$  is now of degree $r_j(n-1)$.  Since 
the differential operators $S_{n,x_j}$ pairwise commute (see Fact $\ref{faitselem}$ $(ii)$), coming back to the definition of $P_l$, one gets
$$t^kP_l(t)=S_{1,x_j}^{(s_j)}\circ S_{1,x_j}^{(r_j-s_j)}\circ V(S_{1,x_j})\circ [t^k]{\circ_{j'\neq j}S_{n,x_{j'}}^{(r_{j'})}\left(A_l(t)\right)\enspace.}$$
But $A_l(t)$ vanishes at $\alpha_i$ at order $\geq {{\rho}} n$ and $S_{1,x_j}^{(r_j-s_j)}\circ V(S_{1,x_j})\circ [t^k]{\circ_{j'\neq j}S_{n,x_{j'}}^{(r_{j'})}}$ involves differentials of order at most ${{\rho}} n-s_j$, hence, by the Leibniz formula, 
$S_{1,x_j}^{(r_j-s_j)}\circ V(S_{1,x_j})\circ [t^k]{\circ_{j'\neq j}S_{n,x_{j'}}^{(r_{j'})}(A_l(t))}$ vanishes at $\alpha_i$ with multiplicity at least $s_j\geq 1$.  Now, by definition, $\varphi_{\alpha_i,x_j,s_j}\circ S_{1,x_j}^{(s_j)}=\varphi_{i,x_j,0}$ whose kernel contains the ideal $(t-\alpha_i)$,  and thus $\varphi_{\alpha_i,x_j,s_j}(t^kP_{l}(t))=0$ which is precisely what is claimed. \qed

{\subsection{Valuations}\label{valu}
We start with a couple of elementary linear algebra  and valuations remarks.
\begin{facts} \label{valufacts} Let $A=\sum_{i\geq k}A_i/z^i\in M_n(K[[1/z]])$ with $k\geq 0$. Then $\det(A)$ is of order ${\rm ord}_{\infty}(\det(A))\geq kn$ and the coefficient of order $kn$ is $\det(A_k)$.
\end{facts}
\begin{lemma}
\label{linalg}
\begin{itemize}
\item[$($i$)$] Let $A\in M_n(K[z])$ a $n\times n$ matrix with polynomial entries. We assume
\begin{itemize} 
\item[$(a)$] The degree of the first row $($maximum of the degrees of the entries$)$ is at most $m\geq 0$.
\item[$(b)$] There exists a lower triangular unipotent element $S\in M_n(K[[1/z]])$ such that the $i$-th rows $L_i$ of $SA$ have an order ${\rm ord}_{\infty}\, L_i\geq k$ for $2\le i \le n$ $($note that the first row $L_1$ of $SA$ is the first row of $A$$)$.
\item[$(c)$] One has $k(n-1)\geq m+1$.
\end{itemize}
Then, $$\det(A)=0\enspace.$$
\item[$($ii$)$] Let $L_i(z)=\sum_{j\geq k}L_i^{(j)}z^{-j}$ the power series expansion of $L_i$ for $i\geq 2$ and $L_1=\sum_{j=0}^{m}L_1^{(j)}z^j$. Assume $(a)$, $(b)$ and $k(n-1)=m$. Then, 
$$\det(A)=\det({}^t\kern-1ptL_1^{(m)},{}^t\kern-1ptL_2^{(k)},\ldots,{}^t\kern-1ptL_n^{(k)})\enspace.$$
\end{itemize}
\end{lemma}

\proofname

Denote by  $D=(D_1,\ldots,D_n)$ the  row vector where $D_t$ is the $(1,t)$-th cofactor of $SA$.  Developing the determinant of $SA$ along the first row yields 
$$\det(SA)=<L_1, D>\enspace,$$
where $<\cdots,\cdots>$ is the canonical bilinear form on $K[z]^n$. Thus, by hypothesis $(a)$ and $(b)$, the order ${\rm ord}_{\infty}(\det(SA))\geq (n-1)k-m\geq 1$ by $(c)$. 
On the other hand, by invariance of the determinant by similarity, $\det(A)=\det(SA)$ is a polynomial, hence if it is non-zero, its order is necessarily $\leq 0$; this yields $(i)$. Point $(ii)$ follows similarly by expanding the first row as a polynomial in $z$ and the other rows as power series in $1/z$.  \qed

\begin{lemma}\label{antinew} 
{Let $A$ be a commutative unitary ring and $\{a_k\}_{k\ge 0}\subset A$ a sequence of elements of $A$. 
Let $f: A[x,y]\longrightarrow A[x]$ be the $A[x]$-linear map defined by $f(y^k)=a_k\in A$ and, for $\lambda\in A$, $g$ the $A[y]$-linear map $A[x,y]\longrightarrow A[y]$ defined by $g(x^k)=\lambda a_k$. 
Let $\tau$ be the transposition inverting $x$ and $y$. Then, for $P\in A[x,y]$ which is antisymmetric $(${\it i.e.} $\tau(P)=-P)$, 
$$f\circ g(P)=0\enspace.$$} 
\end{lemma}

\proofname

It is obvious. \qed

\begin{lemma}\label{racinesymetrique2}
Let $A$ be a commutative unitary ring, $P(X, Y_1,\ldots, Y_{{{\rho}}})\in A[\alpha,X,Y_1,\ldots,Y_{{{\rho}}}]$ be a polynomial with $\tau(P)=-P$ for any transposition $\tau$ interverting $X$ and one of the $Y_i$. Assume moreover that $(X-\alpha)^T \mid P$ for some $T\geq 1$. We further assume that we are given~$:$
 \begin{itemize}
 \item[$(i)$] A differential operator $\partial=c_0+c_1\frac{\partial}{\partial X}+\cdots+c_l\frac{\partial^l}{\partial X^l}$  of order $\leq T-1$ with $c_i\in A$.
 \item[$(ii)$] {{An $A[\alpha,X]$-}}linear map $f: A[\alpha,X,Y_1,\ldots,Y_{{{\rho}}}]\longrightarrow A[\alpha,X,Y_1,\ldots,Y_{{{\rho}}}]$.
 \item[$(iii)$] For each $1\leq i\leq {{\rho}}$, a pair of linear maps $f_i,g_i$  satisfying the hypothesis of Lemma~$\ref{antinew}$ for the pair of variables $(X,Y_i)$.
 \end{itemize} 
 Then, $$\ev_{X\rightarrow\alpha}\circ(f_1\circ g_1+\cdots+f_{{{\rho}}}\circ g_{{{\rho}}}+\partial \circ f)(P)=0\enspace.$$
\end{lemma}

\proofname

Indeed, $\sum_{i=1}^{{{\rho}}}f_i\circ g_i(P)=0$ by Lemma~\ref{antinew} and $f(P)$ is still divisible by $(X-\alpha)^T$ since $f$ treats both $X$ and $\alpha$ as scalars.  Finally, since the order of $\partial$ is $\leq T-1$, $X-\alpha$ divides $\partial(f(P))$, hence $\ev_{X\rightarrow \alpha}\circ\partial\circ f (P)=0$. \qed

The following lemma is a variant of derivation of a parameter integral. We write as a formal lemma that works for more general operators.
\begin{lemma}
\label{dericommutation}
The property holds for the differential operators~$:$
\begin{enumerate}
\item[$(i)$] $\frac{\partial}{\partial \alpha}\circ [\alpha]=\id+[\alpha]\circ \frac{\partial}{\partial \alpha}.$
\item[$(ii)$] $\frac{\partial}{\partial \alpha}\circ\ev_{X_k\rightarrow\alpha}=\ev_{X_k\rightarrow\alpha}\circ \frac{\partial}{\partial \alpha}+\ev_{X_k\rightarrow\alpha}\circ \frac{\partial}{\partial X_k}$.
\item[$(iii)$] $\frac{\partial}{\partial \alpha}\circ\inte_{x,X_k}=\inte_{x, X_k}\circ \frac{\partial}{\partial \alpha}$.
\item[$(iv)$] $S_{1,x,X_k}\circ [1/X_k]=[1/X_k]\circ S_{1,x,X_k}-[1/X_k]$.
\end{enumerate}
\end{lemma}

\proofname

These identities follow from the definitions. \qed

\section{Pad\'e approximation}
\label{deter}
\subsection{Non-vanishing of the crucial determinants}
\label{debutdeter}

For a positive integer $n$ and non-negative integer $l$ with {$0\le l \le {{\rho}} m$}, recall the polynomials  
{$P_{l}(z), P_{l,i,s_j}(z)$} defined in $(\ref{Qnl})$ and $(\ref{Qnlijsj})$ respectively.
We define a column vector $\boldsymbol{p}_{l}(z)\in K[z]^{{{\rho}} m+1}$ by
\begin{align*}
\boldsymbol{p}_{l}(z)={}^t\!\Biggl(P_{l}(z),P_{l,i,s_j}(z)\Biggr)_{\kern-3pt\substack{1\le i \le m ,  \\ 1\le j \le d, 1\le s_j \le r_j}}\enspace.
\end{align*}
The aim of this subsection is to prove the following proposition.
\begin{proposition} \label{non zero det}
\begin{align} \label{original problem}
\Delta(z)=
{\rm{det}}{\left(\rule{0pt}{11pt}
\boldsymbol{p}_{0}(z) ,\ldots , \boldsymbol{p}_{{{\rho}} m}(z)
\right)}
\in K\setminus\{0\}\enspace. 
\end{align}
\end{proposition} {
In the next subsections we give the proof of Proposition \ref{non zero det}, which is quite involved. However, it follows the broad lines of the non-vanishing proof in our previous paper \cite{DHK3} which we refer to as a tem\-plate. 
The previously considered case (a single $x$) being much simpler, a few steps happened to be more straightforward.

{
The first reduction is performed in Subsection $4.1.1$, where we prove that $\Delta(z)\in K$ ({\it i.~e.} $\Delta(z)$ is a constant independent of $z$). 
This crucial step is summed up by Lemma~$\ref{sufficient condition}$. This is achieved by a valuation argument, by proving that $\Delta(z)$ is a polynomial in $z$, but of negative valuation with respect to $z$ (thanks to Lemma~\ref{linalg}). 

Second, we move on to express $\Delta=\Delta(z)$ as the image of a given polynomial by our base operators $\varphi$. This is carried out  in Subsection $4.1.2$: the needed result is described by Lemma~$\ref{another presentation}$.

We then consider factor $\Delta$ viewing the $\alpha_i$ as variables: this is done in Section $4.1.4$, the main result being stated as Proposition~$\ref{decompose Cnum}$.  
To achieve this step, we introduce ``extra'' variables (where the variable $t$ is split in many variables {${{t_{i,s_j}}}$} as there are columns). Basically this follows the principle known since A.~T.~Vandermonde. After showing that the determinant is in fact the determinant of a given minor of the original matrix (all the others canceling out, this is just an extension of the previous valuation argument), we take advantage of the false variables and of the linearity of the operators $\varphi$ to express the determinant as desired.
It can be shown, using Leibniz formula and multilinearity of the determinant, that the few spurious terms which would not be nicely factored as desired, actually cancel out.
After having checked homogeneity and found the trivial monomial factors in $\alpha_i$ (Lemma~$\ref{homo}$), we need to show that $\Delta$ also factors through the $\alpha_i-\alpha_j$ at the appropriate power. 

All these steps, up to minor adjustments are essentially the same as the template case of a single differential shift considered in \cite{DHK3}.

To achieve the final factorization, one shows that the derivative of $\Delta$ with respect to any one of the variables $\alpha_i$ vanishes at the other $\alpha_j$ at the appropriate power. Unfortunately, derivation and operators $\varphi$ do not commute properly. We thus measure the defect of commutativity (Lemma~$\ref{prelimiphi}$), and proceed to note that a derivative of sufficiently high order will look like a derivative, because the operator $\varphi$ is essentially an iterated primitive. 
Since we start with a polynomial vanishing of a high order along $\alpha_i-\alpha_j$, these derivatives tend to produce a lot of vanishing ({{Lemmata~\ref{antinew} and~\ref{racinesymetrique2}}}). However, this argument is not quite enough since derivatives of lower order are not enough to remove all the primitivation, built in the operators $\varphi$. It is taken care via a symmetry argument (Lemma~\ref{antinew}) and this finally shows the corresponding integrals all vanish.

 {Up to this stage, the argument is close to the principle of the proof of \cite{DHK3}. However,  just replicating it would not provide enough multiplicities.  One therefore performs an interpolation between the maps $\inte_{x_j}$ and $\inte_{x_{j'}}$ in order to recover the otherwise undetected multiplicities (Lemma~$\ref{higherintegrals}$). Indeed, it can be shown, again by symmetry, that  the error term coming from the interpolation vanishes.}

After having performed the required interpolation, one is ready to prove factorization: a combinatorics argument is enough to conclude. 
{This is done in two steps: firstly to prepare the conditions of Lemma ~$\ref{condicombi}$ by combinatorial argument, and secondly to 
clearly state the conclusion of Lemma ~$\ref{condicombi}$ (it turns out that once an interpolation is done, } one is actually reduced to the same combinatorial problem
in the single shift case, which has already been laid out in our previous
paper~\cite{DHK3}).

Finally, it will remain to check that the last numerical constant (which depends on $\boldsymbol{x}$ as well as of $r,m$ and $n$) does not vanish. This is spread in two subsections (4.1.4 and 4.1.5). 
Firstly, an induction step (Lemma~$\ref{induction}$) on the number $m$  of algebraic values $\alpha_1,\ldots,\alpha_m$  at which we evaluate the generalized polylogarithmic functions allows to reduce to the case $m=1$. In this more simple situation, the remaining constant is shown to be an Hermite determinant in the variables $x_1,\ldots,x_d$: this is done in a relatively classical way 
  in Lemmata $\ref{Q}$, $\ref{H}$, $\ref{N tilde}$, $\ref{normaconstante}$ and $\ref{last}$, to complete the proof of Proposition~$\ref{non zero det}$. }{
  
  It should be mentioned that the whole construction is effective, and that the actual value of the determinant can be computed if needed.
  }}

\subsubsection{First step} 
{
For a positive integer $n$ and non-negative integer $l$ with $0\le l \le {{\rho}} m$, recall the polynomials 
$P_{l}(z), P_{l,i,s_j}(z)$ defined in $(\ref{Qnl})$ and $(\ref{Qnlijsj})$ respectively.}

\begin{lemma} \label{sufficient condition} 
Define a column vector ${\boldsymbol{u}}_{l}\in K^{{{\rho}} m}$ by
\begin{align*}
{{\boldsymbol{u}}}_{l}={}^t\!\Biggl(\varphi_{\alpha_i,x_j,s_j}(t^nP_{l}(t))\Biggr)_{\substack{1\le i \le m\\ 1\le j \le d, 1\le s_j \le r_j}}\enspace.
\end{align*} 
We use the same notations as in Proposition $\ref{non zero det}$.
Then{, there exists a non-zero element $c\in K$ such that}
$ 
\Delta(z)={c\cdot\det({{\boldsymbol{u}}}_0,\ldots,{{\boldsymbol{u}}}_{{{\rho}} m-1})}\in K. 
$
\end{lemma}

\proofname

Let $$c=\dfrac{1}{({{\rho}} m(n+1))!}\left(\dfrac{d}{dz}\right)^{{{\rho}} m(n+1)}\kern-7ptP_{{{\rho}} m}(z)\enspace,$$ be the coefficient of highest degree ($={{\rho}} mn+{{\rho}} m$) of the polynomial $P_{{{\rho}} m}(z)$. 
Consider
$$\left(\kern-5pt\begin{array}{cccc}
1 & 0 & \cdots& 0\\
\Phi_{1}(x_1,\alpha_1/z)  & -1 & 0 & \vdots\\
\vdots & \ddots & \ddots & \vdots\\
\Phi_{ r_d}(x_d,\alpha_m/z) & 0 & \cdots & -1
\end{array}\kern-5pt\right)\kern-3pt {{\left(\kern-5pt\begin{array}{ccc} \boldsymbol{p}_0(z) ,\ldots, \boldsymbol{p}_{{{\rho}} m}(z)\\	
	\end{array}\kern-5pt\right)}}=\kern-2pt\left(\kern-5pt\begin{array}{ccc} P_0(z)& \ldots & P_{{{\rho}} m }(z)\\
	R_{0,1,1}(z)&\cdots &R_{{{\rho}} m, 1, 1}(z)\\ \vdots & \ddots & \vdots\\ R _{0,m,r_d}(z) & \cdots & R_{{{\rho}} m,m,r_d}(z)\kern-5pt\end{array}\right)\enspace.$$
Since the entries of the first row are (by definition on $P_l(z)$) polynomials of degree $\leq {{\rho}} mn+{{\rho}} m$ and entries of the other rows
(by Theorem~$\ref{Pade appro Lerch}$) are of valuation at least $n+1$, we can apply Lemma~$\ref{linalg}$ $(ii)$.  We need only to check the coefficients of highest degree (for the first row) 
and of minimal valuation {(for all the other rows $1\le s_j \le r_j, 1\le j \le d.$)}. 
By construction, the vector of highest degree ($={{\rho}} mn+{{\rho}} m$) for the first row is $(0,\ldots,0,c)z^{{{\rho}} mn+{{\rho}} m}$  and 
\begin{align*}
R_{l,i,s_j}(z)=\sum_{k=n}^{\infty}\dfrac{{\varphi_{\alpha_i,x_j,s_j}}(t^kP_{l}(t))}{z^{k+1}}\enspace,
\end{align*} for $0\le l \le {{\rho}} m, \ 1\le i \le m, \ 1\le j \le d$ and $1\le s_j \le r_j$. 
So, by Lemma~$\ref{linalg}$ $(ii)$
$$\det(\boldsymbol{p_0}(z),\ldots,\boldsymbol{p}_{{{\rho}} m}(z))=\pm c\cdot \det\left( \varphi_{\alpha_i,x_j,s_j}(t^nP_l(t))\right)_{\substack{0\leq l\leq {{\rho}} m-1\\ 1\leq i\leq m, 1\leq j\leq d, 1\leq s_j\leq r_j}}\enspace,$$
as claimed. \qed

\subsubsection{Second step} 
{We now define} a column vector ${{\boldsymbol{w}}}_{l}\in K^{{{\rho}} m}$ by ({recall the polynomial $A_l$ was defined in equation~($\ref{beta}$)})
\begin{align*}
&{{\boldsymbol{w}}}_{l}={}^t\!\Biggl(\varphi_{\alpha_i,x_j,s_j}(t^nA_{l}(t))\Biggr)_{\substack{1\le i \le m \\ 1\le j \le d, 1\le s_j \le r_j}}\enspace.
\end{align*}
\begin{lemma} \label{another presentation}
We use the notations as above. Then there exists an {effectively computable} element $E\in \Q[x_1,\ldots,x_d]\setminus \{0\}$ satisfying the following equality~$:$
\begin{align} \label{another rep}
\Theta=E \cdot
{\rm{det}}\,  
{{{\begin{pmatrix}
{{\boldsymbol{w}}}_{0} ,\ldots,
{{\boldsymbol{w}}}_{{{\rho}} m-1} 
\end{pmatrix}}}}\enspace.
\end{align} 
\end{lemma}

\proofname

We fix some $x\in\qu$ which is not negative integer. Note that by definition, the linear operator $S_{1,x_j}=S_{1,x}+(x_j-x)\id$ and first look at $t^nP_l(t)$.
By definition, this is ${[t^n]\circ_{j=1}^dS_{n,x_j}^{(r_j)}(A_l(t))}$. 

Setting 
$$H(T)=\frac{1}{n!}\prod_{l=1}^n(T-l)\enspace, \ \ \ H_x(T)=\prod_{k=1}^dH(T+x_k-x)^{r_k}\enspace.$$ 
By Lemma~$\ref{key1}$, we have
$$[t^n]\circ S_{n,x_j}=H(S_{1,x_j})\circ [t^n]\enspace, \ \ \ {[t^n]\circ_{j=1}^dS_{n,x_j}^{(r_j)}=H_x(S_{1,x})\circ [t^n]}\enspace.$$
We apply this remark to a coordinate {$\varphi_{\alpha_i,x_j,s_j}(t^nP_l(t))$} of the vector ${{\boldsymbol{u}}}_l${, ({\it id est} for some  choice of $(i,s_j)$),} and get 
$$\varphi_{\alpha_i,x_j,s_j}(t^nP_l(t))=\varphi_{\alpha_i,x_j,s_j}\circ H_{x_j}(S_{1,x_j})\circ [t^n](A_l(t))\enspace.$$
We now write $H_{x_j}(T)=H_{x_j,s_j}(T)+T^{s_j}Q_{x_j,s_j}(T)$ where {{$H_{x_j,s_j},Q_{x_j,s_j}\in K[T]$ with ${\rm{deg}}\, H_{x_j,s_j}<s_j$}} via the euclidean division. 
Note that $H_{x_j}(0)=H_{x_j,s_j}(0)$.
By definition of $\varphi_{\alpha_i,x_j,s_j}$, we thus have
$$\varphi_{\alpha_i,x_j,s_j}(t^nP_l(t))=\varphi_{\alpha_i,x_j,s_j}\circ H_{x_j,s_j}(S_{1,x_j})\circ [t^n](A_l(t))+ \varphi_{\alpha_i,x_j,0}\circ Q_{x_j}(S_{1,s_j})\circ [t^n](A_l(t))\enspace,$$
by the vanishing Lemma~${\ref{kernel}}$, $\varphi_{\alpha_i,x_j,0}\circ Q_{x_j}(S_{1,s_j})\circ [t^n](A_l(t))=0$ since we differentiate at order\footnote{Indeed, $\deg Q_{x_j}={{\rho}} n-s_j<{{\rho}} n$.} ${< {{\rho}} n}$ a polynomial ($=t^nA_l(t)$) vanishing at order $\geq {{\rho}} n$ at $\alpha_i$. 
Hence,
$$\det({{\boldsymbol{u}}}_l)_{0\leq l\leq {{\rho}} m-1}=\det\left(\varphi_{\alpha_i,x_j,s_j}\circ H_{x_j,s_j}(S_{1,x_j})\circ [t^n](A_l(t))\rule{0mm}{5mm}\right)_{\substack{0\leq l\leq {{\rho}} m-1 \\ 1\leq i\leq m, 1\leq j\leq d, 1\leq s_j\leq r_j}}\enspace.$$
{{By the properties of determinant}}, we get
$$\det({{\boldsymbol{u}}}_l)_{0\leq l\leq {{\rho}} m-1}=\prod_{1\leq j\leq d, 1\leq s_j\leq r_j}H_{x_j,s_j}(0)\cdot \det\left(\varphi_{\alpha_i,x_j,s_j}\circ [t^n](A_l(t))\rule{0mm}{5mm}\right)_{\substack{0\leq l\leq {{\rho}} m-1 \\ 1\leq i\leq m, 1\leq j\leq d, 1\leq s_j\leq r_j}}\enspace.$$
In order to complete the proof, we need only to prove
$$E:=\prod_{1\leq j\leq d, 1\leq s_j\leq r_j}H_{x_j,s_j}(0)^m\in \qu[x_1,\ldots,x_d]\setminus\{0\}\enspace.$$
We thus just need to compute the constant term of the polynomials $H_{x_j,s_j}$. 
Since the constant term of $H_{x_j,s_j}$ is
$$H_{x_j,s_j}(0)=H_{x_j}(0)=\prod_{k=1}^dH(x_k-x_j)^{r_k}\enspace,$$
and $E$ is non-zero since $x_l-x_k\not\in \zu$ for $1\leq k\neq l\leq d$. This concludes the proof of the lemma. \qed

\subsubsection{Third step} 
Now, we consider the ring\footnote{Recall that the index $s_j$ runs through the tagged set of integers {
$\{_j\kern-1pt1,_j\kern-1pt2,\ldots\}
\simeq
\{(j,1),(j,2),\ldots,\}
$}.}}
${K[t_{i,s_j}]}_{\substack{1\leq i\leq m \\ 1\leq j \leq d, 1\leq s_j\leq r_j}}$, the ring of polynomials in ${{\rho}} m$ variables over $K$. 
For each variable {$t_{i,s_j}$}, {any} $\alpha\in K \setminus\{0\}$, {$x_j\in K$} which is not negative integer, one has a well defined map {$\varphi_{\alpha_i,x_j,s_j}=\varphi_{\alpha_i,x_j,t_{i,s_j}}$} for $1\le i \le m$~: 
\begin{equation}
\label{defvarphi}
\varphi_{\alpha_i,x_j,s_j}:K[t_{i',{s_j}'}]_{\substack{(i',s'_j)\neq (i,s_j)}}[t_{i,{s_j}}]\longrightarrow K[t_{i',s'_j}]_{(i',s'_j)\neq (i,s_j)}; \ \ t^k_{i,s_j}\mapsto \dfrac{\alpha_i^{k+1}}{(k+x+1)^{s_j}}\enspace,
\end{equation}
using the definition above where $K[t_{i',{s'_j}}]_{\substack{1\leq i' \leq m\\  1\leq {s'_j}\leq r_j}}$ is seen as the one variable polynomial ring $K'[t_{i,s_j}]$ over $K'=K[t_{i'',s''_j}]_{(i'',s''_j)\neq (i,s_j)}$.
Then we have 
{$${{\boldsymbol{w}}}_l={}^t\!\Biggl(\varphi_{\alpha_i,x_j,s_j}(t^n_{i,s_j}{A}_{l}(t_{i,s_j}))\Biggr)_{\substack{1\le i \le m \\ 1\le j \le d, 1\le s_j \le r_j}}\enspace.$$}
We now define for non-negative integers $u,n$
\begin{align*}
{{\hat{P}}}_{u,n}(t_{i,s_j})={{\hat{P}}}_{u}(t_{i,s_j})
                 :=\prod_{i=1}^{m}\prod_{j=1}^{d}\prod_{s_j=1}^{r_j}\left[t_{i,s_j}^u\prod_{k=1}^m(t_{i,s_j}-\alpha_k)^{{{\rho}} n}\right]{\left(\prod_{(i_1,s_{j_1})<(i_2,s_{j_2})}(t_{i_2,s_{j_2}}-t_{i_1,s_{j_1}})\right)}\enspace,
\end{align*}
where the order $(i_1,s_{j_1})<(i_2,s_{j_2})$ means lexicographical order\footnote{ $(i_1,s_{j_1})<(i_2,s_{j_2})$ if $i_1<i_2$ or $i_1=i_2$ and $j_1<j_2$ or $(i_1,j_1)=(i_2,j_2)$ and $s_{j_1}<s_{j_2}$.} and denote $(t_{i,s_j})$ by $\boldsymbol{t}$. {If $\boldsymbol{k}=(k_{i,s_j})$ is a ${{\rho}} m$-tuple of integers, we denote by $\boldsymbol{t}^{\boldsymbol{k}}$ the product $\prod_{i,s_j}t_{i,s_j}^{k_{i,s_j}}$.}
Also set (when no confusion is deemed to occur, we omit the subscripts $\boldsymbol{\alpha}=(\alpha_1,\ldots,\alpha_m)$, $\boldsymbol{x}=(x_1,\ldots,x_d)$)~:
{\begin{equation}\label{defpsi}\psi=\psi_{\boldsymbol{\alpha},\boldsymbol{x}}={{\circ}}_{i=1}^{m}{{\circ}}_{j=1}^{d}{{\circ}}_{s_j=1}^{r_j}\varphi_{\alpha_i,x_j,s_j}\enspace.\end{equation}}
Note that, by the definition of 
$
{\rm{det}} \, 
{{{\begin{pmatrix}
{{\boldsymbol{w}}}_{0} ,\ldots, {{\boldsymbol{w}}}_{{{\rho}} m-1}
\end{pmatrix}}}},
$
we have
\begin{align} \label{rn}
&{\rm{det}}\,  
{{{\begin{pmatrix}
{{\boldsymbol{w}}}_{0} ,\ldots, {{\boldsymbol{w}}}_{{{\rho}} m-1}
\end{pmatrix}}}}=\psi({{\hat{P}}}_{n,n}(\boldsymbol{t}))\enspace.
\end{align}
{Indeed, since $\varphi_{\alpha_i,x_j,s_j}$ treats all variables $t_{i',s'_j}$  except $t_{i,s_j}$ as scalars,
$$\det({{\boldsymbol{w}}}_{0},\ldots,{{\boldsymbol{w}}}_{{{\rho}} m-1})= \psi\left(\prod_{i,s_j}t_{i,s_j}^n(t_{i,s_j}-\alpha_i)^{{{\rho}} n}\sum_{\sigma\in \mathfrak{S}_{{{\rho}} m}}\varepsilon(\sigma)\boldsymbol{t}^{\sigma({\boldsymbol{l}})}\right)\enspace,$$
and the last sum $\sum_{\sigma\in \mathfrak{S}_{{{\rho}}m}}\varepsilon(\sigma)\boldsymbol{t}^{\sigma({\boldsymbol{l}})}$ is nothing but the Vandermonde determinant in $t_{i,s_j}$.}
In this subsection, we prove 
\begin{proposition} \label{decompose Cnum}
We use above notations. Then there exists a constant $c_{u,m}\in \Q[x_1,\ldots,x_d]$ with
\begin{align} \label{tenkai Cnum}
C_{u,m}:=\psi({{\hat{P}}}_{u}(\boldsymbol{t}))=c_{u,m}\kern-3pt\left(\prod_{i=1}^m\alpha_i\right)^{\kern-5pt {{\rho}}(u+1)+{{\rho}}^2n+\left(\sum_{1\le j_1<j_2\le d}r_{j_1}r_{j_2}+\sum_{j=1}^d{r_j\choose 2}\right)}\kern-15pt\prod_{1\le i_1<i_2\le m}\kern-12pt(\alpha_{i_2}-\alpha_{i_1})^{(2n+1){{\rho}}^2}\enspace.
\end{align}
\end{proposition}
It is also easy to see that since all the variables $t_{i,s_j}$ have been specialized, $C_{u,m}\in K$ is a polynomial in the $\alpha_i,x_j$. The statement is then about a factorization of this polynomial. 
To prove Proposition $\ref{decompose Cnum}$, we are going to perform the following steps~:
\begin{itemize}
\item[$(i)$] Show that $C_{u,m}$ is homogeneous {in $\boldsymbol{\alpha}$} of degree $${{\rho}} m(u+1)+{{\rho}}^2m^2n+{{\rho}}^2\binom{m}{2} 
+m\left(\displaystyle{\sum_{1\le j_1<j_2\le d}}r_{j_1}r_{j_2}+\displaystyle{\sum_{j=1}^d}{r_j\choose 2}\right)\enspace.$$
\item[$(ii)$] Show that $\prod_{i=1}^m\alpha^{{{\rho}} (u+1)+{{\rho}}^2n+\left(\sum_{1\le j_1<j_2\le d}r_{j_1}r_{j_2}+\sum_{j=1}^d{r_j\choose 2}\right)}_i$ divides $C_{u,m}$.
\item[$(iii)$] Show that $\prod_{1\le i_1<i_2\le m}(\alpha_{i_2}-\alpha_{i_1})^{(2n+1){{\rho}}^2}$ divides $C_{u,m}$.
\item[$(iv)$] {Show that the remaining $c_{u,m}$ (now a constant in $\boldsymbol{\alpha}$) is a polynomial with rational coefficients in $\boldsymbol{x}$.}
\end{itemize}
We first prove $(i)$ and $(ii)$ and $(iv)$.
\begin{lemma} \label{homo}
The number $C_{u,m}$ is homogeneous {in $\boldsymbol{\alpha}$} of degree $${{\rho}} m(u+1)+{{\rho}}^2m^2n+{{\rho}}^2\binom{m}{2} 
+m\left(\displaystyle{\sum_{1\le j_1<j_2\le d}}r_{j_1}r_{j_2}+\displaystyle{\sum_{j=1}^d}{r_j\choose 2}\right)\enspace,$$
and is divisible by
$$\left(\prod_{i=1}^m\alpha_i\right)^{{{\rho}}(u+1)+{{\rho}}^2n+\left(\sum_{1\le j_1<j_2\le d}r_{j_1}r_{j_2}+\sum_{j=1}^d{r_j\choose 2}\right)}\enspace.$$
\end{lemma}

\proofname

First the polynomial $P_{u}(\boldsymbol{t})$ is a homogeneous polynomial with respect to the variables {$\alpha_i,\boldsymbol{t}$ }of degree
$${{\rho}} mu+{{\rho}}^2m^2n+{{\rho}}^2\binom{m}{2} 
+m\left(\displaystyle{\sum_{1\le j_1<j_2\le d}}r_{j_1}r_{j_2}+\displaystyle{\sum_{j=1}^d}{r_j\choose 2}\right)\enspace.$$
{Thus, by} the definition of $\psi$ {(which is degree increasing by ${{\rho}} m$)}, it is easy to see that $C_{u,m}=\psi({{\hat{P}}}_{u}(\boldsymbol{t}))$ is a homogeneous polynomial with respect to the variables $\alpha_i$ of degree 
$${{\rho}} m(u+1)+{{\rho}}^2m^2n+{{\rho}}^2\binom{m}{2}
+m\left(\displaystyle{\sum_{1\le j_1<j_2\le d}}r_{j_1}r_{j_2}+\displaystyle{\sum_{j=1}^d}{r_j\choose 2}\right)\enspace.$$
Second we show the later assertion. By linear algebra{,} {$\varphi_{\alpha_i,x_j,s_j}(Q(t))=\alpha_i \ev_{\alpha_i\rightarrow 1}\circ\varphi_{1,x_j,s_j}(Q(\alpha_i t))$} ({\it  i.e.} the variable $t$ specializes in $1$, {\it confer} Lemma~$\ref{prelimiphi}$ $(ii)$) for  any polynomial $Q(t)\in K[t]$. 
So, by composition, the same holds for $\psi$.
{So, putting $\boldsymbol{1}=(1,\ldots,1)\in K^m$, 
We now compute
$$
{{\hat{P}}}_{u}({{(\alpha_i t_{i,s_j})_{i,s_j}}})=
\left(\prod_{i=1}^m\alpha_i\right)^{{{\rho}}(u+1)+{{\rho}}^2n+\left(\sum_{1\le j_1<j_2\le d}r_{j_1}r_{j_2}+\sum_{j=1}^d{r_j\choose 2}\right)} \cdot  Q_u({\boldsymbol{t}})\enspace,
$$
where 
{\begin{equation}\label{defqu}
\begin{array}{lcl}\displaystyle
Q_u(\boldsymbol{t}) = Q_{n,u,m}(\boldsymbol{t})&= & \displaystyle\left(\prod_{i=1}^m\prod_{j=1}^d\prod_{s_j=1}^{r_j}\left [t_{i,s_j}^u\prod_{k\neq i}(\alpha_i t_{i,s_j}-\alpha_k)^{{{\rho}} n}(t_{i,s_j}-1)^{{{\rho}} n}\right]\right)\\
&\displaystyle \cdot& \displaystyle\rule{0mm}{6mm} \prod_{1\le i_1<i_2\le m}\prod_{(j_1,s_{j_1})<(j_2,s_{j_2})}(\alpha_{i_2}t_{i_2,s_{j_2}}-\alpha_{i_1}t_{i_1,s_{j_1}})\\
& \cdot& \rule{0mm}{6mm} \displaystyle{\prod_{i=1}^m\prod_{(j_1,s_{j_1})<(j_2,s_{j_2})}}\left(t_{i,s_{j_2}}-t_{i,s_{j_1}}\right)\enspace,
\end{array}
\end{equation}}
by linearity, 
\begin{equation} \label{equality 1}
C_{u,m}=\left(\prod_{i=1}^m\alpha_i\right)^{{{\rho}}(u+1)+{{\rho}}^2n+\left(\sum_{1\le j_1<j_2\le d}r_{j_1}r_{j_2}+\sum_{j=1}^d{r_j\choose 2}\right)}\cdot\psi_{\boldsymbol{1},\boldsymbol{x}}(Q_u(\boldsymbol{t}))\enspace.
\end{equation}
This concludes the proof of the lemma.} \qed

Note, by definition we have $C_{u,m}\in \Q[\alpha_1,\ldots,\alpha_m,x_1,\ldots,x_d]$ and thanks to Lemma $\ref{homo}$, we have $c_{u,m}\in \Q[x_1,\ldots,x_d]$ assuming step $(iii)$.

\vspace{1.5\baselineskip}
Now we consider $({{iii}})$. 
The case of $m=1$ is trivial. Thus we assume $m\ge 2$. 

We need to show that $(\alpha_{i_2}-\alpha_{i_1})^{(2n+1){{\rho}}^2}$ {divides} $C_{u,m}$. Without loss of generality, after renumbering, we can assume that $i_2=2,i_1=1$.
To ease notations, we are going to take advantage of the fact that $m\geq2$, and set {$X_{s_j}=t_{1,s_j}, Y_{s_j}=t_{2,s_j}$ }and $\alpha_1=\alpha$, $\alpha_2=\beta$. {We want to treat all the variables $t_{i,s_j}$, for $i\geq 3$ as constants as they will play no subsequent role.}
So our polynomial ${{\hat{P}}}_u$ rewrites as
{\begin{align*}
{{\hat{P}}}_{u}(\boldsymbol{X},\boldsymbol{Y})&=\prod_{j=1}^d\prod_{s_j=1}^{r_j}\left(X_{s_j}Y_{s_j}\right)^u\left[(X_{s_j}-\alpha)(X_{s_j}-\beta)(Y_{s_j}-\alpha)(Y_{s_j}-\beta)\right]^{{{\rho}} n}\\
& \cdot \prod_{(j_1,s_{j_1})<(j_2,s_{j_2})}(X_{s_{j_2}}-X_{s_{j_1}})\prod_{(j_1,s_{j_1})<(j_2,s_{j_2})} \left(Y_{s_{j_2}}-Y_{s_{j_1}}\right)\\
&\cdot \prod_{\substack{(j_2,s_{j_2}) \\ (j_1,s_{j_1})}}(Y_{s_{j_2}}-X_{s_{j_1}})
\cdot {c(t_{i,s_j})}\prod_{k\geq 3}\prod_{j=1}^d\prod_{s_j=1}^{r_j}\prod_{\substack{(j_1,s_{j_1}) \\ (j_2,s_{j_2})}}(t_{k,s_j}-X_{s_{j_1}})(t_{k,s_j}-Y_{s_{j_2}})\enspace,
\end{align*}}
{where $c(t_{i,s_j})$ is shortened for $c(t_{i,s_j})_{i\ge 3}$ and is defined by}
{\begin{align*}
c(t_{i,s_j}):&=\prod_{i=3}^{m}\prod_{j=1}^d\prod_{s_j=1}^{r_j}\left[t_{i,s_j}^u\prod_{k=1}^m(t_{i,s_j}-\alpha_k)^{{{\rho}} n}\right]\cdot \prod_{\substack{i_1,i_2\geq 3 \\ (i_1,  s_{j_1})<(i_2, s_{j_2})} }(t_{i_2,s_{j_2}}-t_{i_1,s_{j_1}})\enspace.
\end{align*}} 
{In order to differentiate computations between the set of variables associated to $\alpha$ and $\beta$ respectively,  we set 
\begin{align*}
&\psi_{\alpha}=\psi_{\alpha,\boldsymbol{x}}={{\circ}}_{j=1}^d{{\circ}}_{s_j=1}^{r_j}\varphi_{\alpha,x_j,s_j}\enspace, \ \ 
\psi_{\beta}=\psi_{\beta,\boldsymbol{x}}={{\circ}}_{j=1}^d{{\circ}}_{s_j=1}^{r_j}\varphi_{\beta,x_j,s_j}\enspace.
\end{align*}}
{In other words, we have
$$\psi=\psi_{\alpha,\boldsymbol{x}}\circ \psi_{\beta,\boldsymbol{x}}\circ \Xi\enspace,$$ where 
{$\Xi={\displaystyle{{{\circ}}_{i\ge3}}}{{\circ}}_{j=1}^d{{\circ}}_{s_j=1}^{r_j}\varphi_{\alpha_i,x_j,s_j}={{\circ}}_{i\geq 3}\psi_{\alpha_i,\boldsymbol{x}}$.}}

We now write ${{\hat{P}}}_{u}$ as a product of a symmetric and an antisymmetric polynomial:
{\begin{align*}
{{f_{u,n}(\alpha,\beta,\boldsymbol{X},\boldsymbol{Y})}}=f(\alpha,\beta,\boldsymbol{X},\boldsymbol{Y}) &:=\prod_{j=1}^d\prod_{s_j=1}^{r_j}\left(X_{s_j}Y_{s_j}\right)^u\left[(X_{s_j}-\alpha)(X_{s_j}-\beta)(Y_{s_j}-\alpha)(Y_{s_j}-\beta)\right]^{{{\rho}} n}\\
&\cdot {c(t_{i,s_j})}\prod_{k\geq 3}\prod_{j=1}^d\prod_{s_j=1}^{r_j}\prod_{\substack{(j_1,s_{j_1}) \\ (j_2,s_{j_2})}}(t_{k,s_j}-X_{s_{j_1}})(t_{k,s_j}-Y_{s_{j_2}})\enspace,
\end{align*}}
{where $c(t_{i,s_j})$ is shortened for $c(t_{i,s_j})_{i\ge 3}$ and}
{\begin{align*}
g(\boldsymbol{X},\boldsymbol{Y})&=g(\alpha,\beta,\boldsymbol{X},\boldsymbol{Y}):=\kern-3pt\prod_{(j_1,s_{j_1})<(j_2,s_{j_2})}(X_{s_{j_2}}-X_{s_{j_1}})\prod_{(j_1,s_{j_1})<(j_2,s_{j_2})} (Y_{s_{j_2}}-Y_{s_{j_1}})\cdot \prod_{\substack{(j_2,s_{j_2}) \\ (j_1,s_{j_1})}}(Y_{s_{j_2}}-X_{s_{j_1}})
\enspace.
\end{align*}}
So that ${{\hat{P}}}_{u}={{\hat{P}}}=fg$ (for the rest of the proof, the indexes $u,n$ will not play any role and may be conveniently left off to ease reading).

{In order to tackle the different shifts, we need to be able to interpolate between them}. 
We have the following elementary observations.

Let $(j,s_j)$ as above, {\it i.e.} $1\leq j\leq d, 1\leq s_j\leq r_j$; we set $k=\sum_{i=1}^{j-1}r_i+s_j$ (as usual, the empty sum is equal to zero), {so that we have established a natural bijection between $\{(j,s_j)\mid 1\leq j\leq d, 1\leq s_j\leq r_j\}$ and $\{1,\ldots,r\}$. }
\begin{lemma}
\label{higherintegrals} 
\begin{itemize}
\item[$({{i}})$] Let $x,y\in\qu$ {{which are not negative integers}} and $n\in \N$. Then we have $$\inte_y=\sum_{k=0}^{n-1}(x-y)^{k}\inte^{(k+1)}_{x}+(x-y)^{n}\inte^{(n)}_x\circ \inte_{y}\enspace.$$

\item[$({{ii}})$] For each $1\le j \le d$, {and any integer $k=\sum_{i=1}^{j-1}r_i+s_j$} {{with $1\le s_j\le r_j$}}, we have
\begin{align*}
{\rm{Prim}}_{x_j}&=\left({{\circ}}_{i=1}^{j-1}\left[(x_i-x_j)^{r_i}\inte_{x_i}^{(r_i)}\right]\circ\inte_{x_j}\right)\\
&+\left(\sum_{{\mu}\leq \sum_{i=1}^{j-1}r_i}\prod_{i=1}^{n-1}(x_i-x_j)^{r_i}(x_n-x_j)^{{{s_n-1}}}{{\circ}}_{i=1}^{n-1}\inte_{x_i}^{(r_i)}\circ \inte_{x_n}^{(s_n)}\right)\enspace,
\end{align*}
where at each stage, one writes ${\mu}=\sum_{i=1}^{n-1}r_i+s_n$ in a unique fashion for some $s_n$ between $1$ and $r_n$, where the primitive is always taken with respect to {the same fixed variable} $X_{s_j}$ $($respectively $Y_{s_j}$ for specialization at $\beta)$. As usual, empty product is $1$, empty sum is zero and empty composition is identity.

In particular, we have  
\begin{equation}
\label{profondeur}
\begin{array}{lcl}\displaystyle
&\varphi_{\alpha,x_j,s_j} \circ\deri_{x_j}^{(s_j-1)} ={\displaystyle [\alpha]\circ\ev_{{{X_{s_j}\rightarrow\alpha}}} \circ_{i=1}^{j-1}\left[(x_i-x_j)^{r_i}\inte_{x_i}^{(r_i)}\right]\circ\inte_{x_j}}\\
&{+ \displaystyle [\alpha]\circ\ev_{{{X_{s_j}\rightarrow\alpha}}} \left(\sum_{\mu \leq \sum_{i=1}^{j-1}r_i}\prod_{i=1}^{n-1}(x_i-x_j)^{r_i}(x_n-x_j)^{{{s_n-1}}}\circ_{i=1}^{n-1}\inte_{x_i}^{(r_i)}\circ \inte_{x_n}^{(s_n)}\right)\enspace,}
\end{array}
\end{equation}
\end{itemize}
\end{lemma}

\proofname

The first property follow from the definition of the map $\inte_x$, by checking on the basis {of the polynomial ring}. 
We now prove $({{ii}})$.  We proceed by induction on $j$.  If $j=1$, there is nothing to prove since {both compositions and sums are empty.}
Thus we may assume  $j\geq 2$ and $(ii)$ to be true for all $1\leq {l} \leq j-1$. 
We now write, using Lemma~$\ref{higherintegrals}$ $(i)$, {with $y=x_j$,  $x=x_1$ and $n=r_1$,}
$${
\inte_{x_j}=\sum_{l=0}^{r_{1}-1}(x_{1}-x_{j})^l\inte_{x_{1}}^{(l+1)}+(x_{1}-x_{j})^{r_{1}}\inte_{x_{1}}^{(r_{1})}\circ\inte_{x_j}\enspace.}
$$
We now apply the induction hypothesis with the set $\{x_2,\ldots, x_j\}$ of length $j-1$ and get~:
{\small{$$
\inte_{x_j}={{\circ}}_{i=2}^{j-1}\left[(x_i-x_{j})^{r_i}\inte_{x_i}^{(r_i)}\right]\circ\inte_{x_j}+\sum_{\mu\leq \sum_{i=2}^{j-1}r_i}\prod_{i=2}^{n-1}(x_i-x_j)^{r_i}(x_n-x_j)^{{s_n-1}}{{\circ}}_{i=1}^{n-1}\inte_{x_i}^{(r_i)} \circ\inte_{x_n}^{(s_n)}\enspace.
$$}}
Plugging in the above relation in the previous relation completes induction and the proof of Lemma~$\ref{higherintegrals}$. \qed

For $\gamma \in K[\alpha,\beta,\underline{X},\underline{Y}]$, we define the substitution morphism
$$\Delta_{\gamma}:K[\alpha,\beta,\underline{X},\underline{Y}]\longrightarrow K[\alpha,\beta,\underline{X},\underline{Y}]; \ \Delta_{\gamma}(\beta)=\gamma\enspace,$$
with $\left. \Delta_{\gamma}\right|_{K[\alpha,\underline{X},\underline{Y}]}={\rm{Id}}_{K[\alpha,\underline{X},\underline{Y}]}$.
\vspace{\baselineskip}
We now rewrite our polynomial.
Let $(j,s_j)$ as above, {\it i.e.} $1\leq j\leq d, 1\leq s_j\leq r_j$; we set $k=\sum_{i=1}^{j-1}r_i+s_j$ (as usual, the empty sum is equal to zero).  We then rewrite our shifts setting $y_k=x_j$ (in other words, each ${{x_j}}$ is repeated ${{r_j}}$ times). Similarly, we rewrite the variable {$X_{s_j}$} as $X_{k}$. 
{In contrast to the first part of the proof which consisted mainly along the same line as Hermite in {seperating} the variable{s}, we shall now need to let both the shifts and the number of primitive to vary with more flexibility. This forces to  take into account: the number of primitivations that are performed, the shifts associated to them and the variable with respect to which primitivation is performed.}

We put for any integer ${{\nu}}$, and any integer $1\leq {{\lambda}} \leq {{\rho}}$, {$$\theta_{\alpha,{{\lambda}},{{\nu}},X_{k} }=\varphi_{\alpha,x_1,{{\nu}}-{{\lambda}},X_{k}}\circ_{i=1}^{{{\lambda}}}\inte_{y_i, X_{k}}\enspace,$$} with as usual empty composition being the identity map, we define in the same way the map theta with specialization at $\beta$ and integration is taken with respect to the variable ${Y_{k}={Y_{s_j}}}$.

We set $\boldsymbol{{{\lambda}}}=({{\lambda}}_{1},\ldots,{{\lambda}}_{{{\rho}}})$, $\boldsymbol{{{\nu}}}=({{\nu}}_{1},\ldots,{{\nu}}_{{{\rho}}})\in \Z^{{{\rho}}}$ and denote by 
\begin{equation}
\Theta_{\alpha,\boldsymbol{{{\lambda}}},\boldsymbol{{{\nu}}}}=\circ_{k=1}^r\theta_{\alpha,{{\lambda}}_{k},{{\nu}}_{k},X_{k}}\enspace, 
\label{defTheta}
\end{equation} and in a same fashion for specialization at $\beta$.

\begin{lemma} \label{Cu2}
We have {$C_{u,m}=C({\boldsymbol{y}})\Xi\circ\Theta_{\beta,\boldsymbol{l},\boldsymbol{l}}\circ\Theta_{\alpha,\boldsymbol{l},\boldsymbol{l}}(fg)$}, where $C(\boldsymbol{y})$ is some non-zero real number depending only on 
$\boldsymbol{y}$, where $\boldsymbol{y}=(y_1,\ldots,y_{{{\rho}}})$ and $\boldsymbol{l}=(1,2,\ldots,{{\rho}})\in \Z^{{{\rho}}}$.
\label{reductheta}
\end{lemma}

\proofname

We {recall} $\boldsymbol{k}:=(1,\ldots,r_1,\ldots,1,\ldots,r_d)\in \Z^{{{\rho}}}$ {and that 
$C_{u,{{m}}}=\Xi\circ \psi_{\alpha,\boldsymbol{x}}\circ\psi_{\beta,\boldsymbol{x}}(fg).$}
It is thus enough to show {$$\psi_{\alpha,\boldsymbol{x}}\circ \psi_{\beta,\boldsymbol{x}}(fg)={C({\boldsymbol{y})}}\Theta_{\alpha,\boldsymbol{l},\boldsymbol{l}}\circ \Theta_{\beta,\boldsymbol{l},\boldsymbol{l}}(fg) \enspace.$$ }
Let $k$ be an integer with $1\leq k\leq {{\rho}}$ and write $k=\sum_{i=1}^{j-1}r_i+s_j$ for $1\leq s_j\leq r_j$ (possible in a unique way). 
We now recall 
$${{\varphi_{\alpha,x_j,s_j,X_{s_j}}=[\alpha]\circ {\rm{Eval}}_{{{X_{s_j}\rightarrow \alpha}}}\circ \inte_{x_j}\circ\inte_{x_j}^{(s_j-1)}}}\enspace.$$ 
We then apply Lemma~$\ref{higherintegrals}$, part $(ii)$, and get 
{\begin{align*}
&\varphi_{\alpha,x_j,s_j,X_{s_j}}=
[\alpha]\circ\ev_{{{X_{s_j}\rightarrow\alpha}}}\kern-3pt\\
&{\left(\kern-3pt\circ_{i=1}^{j-1}\left[(x_i-x_j)^{r_i}\inte_{x_i}^{(r_i)}\right]\kern-3pt\circ\inte_{x_j}\kern-1pt+\kern-3pt\sum_{\mu \leq \sum_{i=1}^{j-1}r_i}A_\mu\circ_{i=1}^{n}\inte_{x_i}^{(r_i)}\circ \inte_{x_n}^{(s_n)}\right)\circ
\inte_{x_j}^{(s_j-1)}\enspace,}
\end{align*} 
for some $A_{\mu}\in\qu[\boldsymbol{y}]\setminus\{0\}$.}
Plugging in the definition of $\theta$ operators, one gets 
\begin{align} \label{varphi tenkai}
\varphi_{\alpha,x_j,s_j,X_{s_j}}=\prod_{i=1}^{j-1}(x_i-x_j)^{r_i}\theta_{\alpha,k,k,X_{k}}+\sum_{\mu\leq \sum_{i=1}^{j-1}r_i}A_\mu\theta_{\alpha,\mu,\mu,X_{k}}\circ\inte_{x_j}^{(s_j-1)}\enspace.
\end{align}
We now claim $\psi_{\alpha,\boldsymbol{x}}(fg)=C(\boldsymbol{y})\Theta_{\alpha,\boldsymbol{l},\boldsymbol{l}}(fg)$. 
Indeed, this can be proven easily by induction on ${{\rho}}$. If ${{\rho}}=1$, there is nothing to prove since 
$\Theta_{\alpha,1,1}=\psi_{\alpha,x_1,1}$. We now assume $\sum_{i=1}^{{j-1}} r_i+s_j<{{\rho}}$ and
$${{{{\circ}}_{l=1}^{s_j}\varphi_{\alpha,x_j,l,X_{l}}}}{{\circ}}_{1\leq  {{i \leq j-1}}}{{\circ}}_{1\leq s_i\leq {{r_i}}}\varphi_{\alpha,x_i,s_i,X_{s_i}}(fg)=C(\boldsymbol{y}){{\circ}}_{\mu=1}^{\sum_{i=1}^{{j-1}} r_i+s_j}\theta_{\alpha,\mu,\mu,X_{\mu}}(fg)\enspace,$$
and prove the same relation for $\sum_{i=1}^{{j-1}} r_i+s_j+1$. 
We separate the computation into two cases,
firstly, $s_j=r_j$. 
In this case, by {{$(\ref{varphi tenkai})$}},
$$\varphi_{\alpha,x_{j+1},_{j+1}\kern-1pt1,X_{_{j+1}\kern-1pt1}}=\prod_{i=1}^j(x_i-x_{j+1})\theta_{\alpha,k+1,k+1,X_{k+1}}+\sum_{\mu\leq \sum_{i=1}^{j}r_i}A_\mu\theta_{\alpha,\mu,\mu,X_{k+1}}\enspace.$$ 
By Lemma~$\ref{antinew}$, $\theta_{\alpha,\mu, \mu,X_{k+1}}\circ\theta_{\alpha,\mu,\mu,X_{m}}(fg)=0$, so, the induction is completed by expanding the composition.
We now turn to the case where $s_j+1\leq r_j$, and we look at $\theta_{\alpha,\mu,\mu,X_{k+1}}\circ\inte_{x_j}^{(s_j)}$ for some $\mu\leq\sum_{i=1}^{j-1}r_i$. Since $\inte_{x_j}=\inte_{y_{\mu+1}}+(y_{\mu+1}-x_j)\inte_{y_{\mu+1}}\circ\inte_{x_j}$, 
and
$$\theta_{\alpha,\mu,\mu,X_{k+1}}\circ\inte_{x_j}^{(s_j)}=\theta_{\alpha,\mu+1,\mu+1,X_{k+1}}\circ\inte_{x_j}^{(s_j-1)}
+(y_{\mu+1}-x_j)\theta_{\alpha, \mu+1,\mu+1,X_{k+1}}\circ\inte_{x_j}^{(s_j)}\enspace,$$
by induction, $\theta_{\alpha,\mu,\mu,X_{k+1}}\circ\inte_{x_j}^{(s_j)}$ is a linear combination of 
$\theta_{\alpha,{{\nu}},{{\nu}},X_{{k}+1}}$ for ${{\nu}}$ varying between $\mu+1$ and $k$, but by Lemma $\ref{antinew}$,
$\theta_{\alpha,{{\nu}},{{\nu}},X_{k+1}}\circ \theta_{\alpha,{{\nu}},{{\nu}},X_{{{\nu}}}}(fg)=0$, hence the induction is also complete in this case and Lemma~$\ref{reductheta}$ is proved. \qed

\begin{lemma} \label{derivative Theta}
Let $\boldsymbol{{{\lambda}}}=({{\lambda}}_{k})_{1\le k \le {{\rho}}},\boldsymbol{{{\nu}}}=({{\nu}}_{k})_{1\le k \le {{\rho}}}$ be two integral vectors in $\zu^{{{\rho}}}$, with $1\leq {{\lambda}}_{k}\leq k$ and $1\leq k \leq {{\rho}}$.
Then 
$$\frac{\partial}{\partial \alpha}\circ\Theta_{\alpha,{\boldsymbol{{\lambda}}},\boldsymbol{{{\nu}}}}=\Theta_{\alpha,{\boldsymbol{{\lambda}}},\boldsymbol{{{\nu}}}}\circ\left( \frac{\partial}{\partial \alpha}-\frac{ry_1}{\alpha}\right)+\sum_{j=1}^{{{\rho}}}\Theta_{\alpha,{\boldsymbol{{\lambda}}},\boldsymbol{{{\nu}}}-\bold{e}_{j}}\circ\left[\dfrac{1}{\alpha}\right]\enspace,$$
where $\bold{e}_{j}:=(0\ldots,0, 1 ,0 \ldots, 0)\in \Z^{{{\rho}}}$ and $1$ is in the $j$-th spot.
\end{lemma}

\proofname

{By definition, $\Theta_{\alpha,{\boldsymbol{{\lambda}}},\boldsymbol{{{\nu}}}}={{\circ}}_{k=1}^{{{\rho}}}\theta_{\alpha,{{\lambda}}_{k},{{\nu}}_{k},X_{k}}$ and  $$\theta_{\alpha,{{\lambda}},{{\nu}},X_{k} }=[\alpha]\circ \ev_{{{X_k\rightarrow \alpha}}}\circ{\inte^{{{({{\nu}}-{{\lambda}})}}}_{x_1,X_k}
\circ_{i=1}^{{{\lambda}}}\inte_{y_i,X_{k}}\enspace.}$$
Thus, using the elementary Lemma~$\ref{dericommutation}$ above,
{\small{$$\renewcommand{\arraycolsep}{3pt}\begin{array}{lcl}\displaystyle
\tfrac{\partial}{\partial \alpha}\circ\theta_{\alpha,{{\lambda}},{{\nu}},X_{k} } & = & {\displaystyle [\tfrac{1}{\alpha}]\circ \theta_{\alpha,{{\lambda}},{{\nu}},X_{k} } +[\alpha]\circ \tfrac{\partial}{\partial \alpha}\circ  \ev_{{{X_k\rightarrow \alpha}}}\circ\inte^{{{({{\nu}}-{{\lambda}})}}}_{x_1, X_k}
\circ_{i=1}^{{{\lambda}}}\inte_{y_i,X_{k}}}\\ 
& = &\rule{0mm}{5mm} \displaystyle [\tfrac{1}{\alpha}]\circ \theta_{\alpha,{{\lambda}},{{\nu}},X_{k} }  +\theta_{\alpha,{{\lambda}},{{\nu}},X_{k} } \circ\tfrac{\partial}{\partial \alpha}+[\alpha]\circ\ev_{X_k\rightarrow \alpha}\circ \tfrac{\partial}{\partial X_k}\circ \inte^{{{({{\nu}}-{{\lambda}})}}}_{x_1, X_k} {\circ_{i=1}^{{{\lambda}}}\inte_{y_i, X_{k}}\enspace.}
\end{array}
$$}}
Since $S_{1,x_1,X_k}=\frac{\partial}{\partial X_k}\circ[X_k]+x_1\id$,
this yields,
{\small{$$\renewcommand{\arraycolsep}{1.2pt}
\begin{array}{lcl}\displaystyle
\kern-10pt\tfrac{\partial}{\partial \alpha}\circ\theta_{\alpha,{{\lambda}},{{\nu}},X_{k} } & = & \displaystyle 
[\tfrac{1}{\alpha}]\circ \theta_{\alpha,{{\lambda}},{{\nu}},X_{k} }  +\theta_{\alpha,{{\lambda}},{{\nu}},X_{k} } \circ\tfrac{\partial}{\partial \alpha}\\ &+ & \displaystyle\rule{0mm}{6mm}[\alpha]\circ\ev_{X_k\rightarrow \alpha}\circ \left(S_{1,x_1,X_k}\circ[\tfrac{1}{X_k}]-[\tfrac{x_1}{X_k}]\right)\circ {\inte^{{{({{\nu}}-{{\lambda}})}}}_{x_1, X_k}
\circ_{i=1}^{{{\lambda}}}\inte_{y_i, X_{k}}}
\\ & = & \displaystyle\rule{0mm}{6mm}[\tfrac{1}{\alpha}]\circ \theta_{\alpha,{{\lambda}},{{\nu}},X_{k} }  +\theta_{\alpha,{{\lambda}},{{\nu}},X_{k} } \circ\tfrac{\partial}{\partial \alpha}
+  \displaystyle\rule{0mm}{5mm}[\alpha]\circ\ev_{X_k\rightarrow \alpha}\circ S_{1,x_1,X_k}\circ[\tfrac{1}{X_k}]\circ {\inte^{{{({{\nu}}-{{\lambda}})}}}_{x_1, X_k}
\circ_{i=1}^{{{\lambda}}}\inte_{y_i, X_{k}}}\\ & - & \rule{0mm}{6mm}\displaystyle [x_1\alpha]\circ\ev_{X_k\rightarrow \alpha}\circ[\tfrac{1}{X_k}]\circ {
\inte^{{{({{\nu}}-{{\lambda}})}}}_{x_1, X_k}
\circ_{i=1}^{{{\lambda}}}\inte_{y_i, X_{k}}}\\
 & = & \rule{0mm}{6mm}\displaystyle\displaystyle\rule{0mm}{6mm}[\tfrac{1}{\alpha}]\circ \theta_{\alpha,{{\lambda}},{{\nu}},X_{k} }  +\theta_{\alpha,{{\lambda}},{{\nu}},X_{k} } \circ\tfrac{\partial}{\partial \alpha}+  \displaystyle\rule{0mm}{5mm}[\alpha]\circ\ev_{X_k\rightarrow \alpha}\circ [\tfrac{1}{X_k}]\circ S_{1,x_1,X_k}\circ \inte^{{{({{\nu}}-{{\lambda}})}}}_{x_1, X_k}
{\circ_{i=1}^{{{\lambda}}}\inte_{y_i, X_{k}}}\\ & - & \rule{0mm}{6mm}\displaystyle [\tfrac{x_1}{\alpha}]\circ\theta_{\alpha,{{\lambda}},{{\nu}},X_k} 
-
[\alpha]\circ\ev_{X_k\rightarrow \alpha}\circ [\tfrac{1}{X_k}]\circ {\inte^{{{({{\nu}}-{{\lambda}})}}}_{x_1, X_k}
\circ_{i=1}^{{{\lambda}}}\inte_{y_i, X_{k}}}\\ 
& 
= & \rule{0mm}{5mm} \rule{0mm}{6mm}\displaystyle\displaystyle\rule{0mm}{5mm}[\tfrac{1}{\alpha}]\circ \theta_{\alpha,{{\lambda}},{{\nu}},X_{k} }  +\theta_{\alpha,{{\lambda}},{{\nu}},X_k } \circ\tfrac{\partial}{\partial \alpha}+  \displaystyle\rule{0mm}{5mm}
[\tfrac{1}{\alpha}]\circ\theta_{\alpha,{{\lambda}},{{\nu}}-1,X_k} -  \displaystyle [\tfrac{x_1}{\alpha}]\circ\theta_{\alpha,{{\lambda}},{{\nu}},X_k} 
-
[\tfrac{1}{\alpha}]\theta_{\alpha,{{\lambda}},{{\nu}},X_k}\\
& = & \displaystyle\rule{0mm}{6mm}
\theta_{\alpha,{{\lambda}},{{\nu}},X_k } \circ\left(\tfrac{\partial}{\partial \alpha}-\tfrac{x_1}{\alpha}\right)+  \displaystyle\rule{0mm}{6mm}\theta_{\alpha,{{\lambda}},{{\nu}}-1,X_k}\circ [\tfrac{1}{\alpha}]\enspace.
\end{array}$$}}
By inputting this in the definition of $\Theta_{\alpha,\boldsymbol{{{\lambda}}},\boldsymbol{{{\nu}}}}$, one obtains the lemma.} \qed

\begin{lemma}\label{racinesymatrique3} 
Put $\boldsymbol{l}=(1,\ldots,{{\rho}})\in \Z^{{{\rho}}}$.
Let ${{\lambda}}, l, k\in\Z$ with $1\le {{\lambda}}, k \le {{\rho}}$ and ${{\lambda}} \le l$. Let $P\in K[{{\alpha,\beta,}}X_{1},\ldots,X_{{{\rho}}},Y_{1},\ldots,Y_{{{\rho}}}]$ be a polynomial such that $(X_{k}-\alpha)^{T_1}(X_{k}-\beta)^{T_2}\mid P$ for non-negative integers $T_1,T_2$ with either $T_1$ or $T_2$ is greater than $1$ and $0\leq l-{{\lambda}}\leq T_1+T_2-1$.
{{Assume the polynomial $P$ satisfies}} $\tau_k(P)=-P$ for transposition $\tau_k$ which {{inverts}} $X_{k}$ and $Y_{k}$. 
Then we have $$\Delta_{\alpha}\circ\theta_{\alpha,{{\lambda}},{{\lambda}}-l,X_{k}}\circ \Theta_{\beta,\boldsymbol{l},\boldsymbol{l}}(P)=0\enspace.$$
\end{lemma}

\proofname

{{Let $({{\eta_i}})_{1\le i \le l}$ be a sequence of elements of $K$ with ${{\eta}}_1=y_{{\lambda}}, {{\eta}}_2=y_{{{\lambda}}-1},\ldots, {{\eta}}_{{\lambda}}=\cdots={{\eta}}_l=y_1$. 
Then there exists a sequence of elements $(b_{h})_{0\le h \le l}$ with $$S^{(l)}_{1,x_1,X_k}=\sum_{h=0}^lb_h\circ_{i=1}^h S_{1,{{\eta}}_i,X_{k}}\enspace,$$
where we mean $\circ_{i=1}^h S_{1,{{\eta}}_i,X_{k}}={\rm{id}}$ if $h=0$. Thus we have
\begin{align}
\theta_{\alpha,{{\lambda}},{{\lambda}}-l,X_{k}}&{=\varphi_{\alpha,x_1,-l,X_{k}}\circ_{u=1}^{{\lambda}}{\rm{Prim}}_{y_u,X_{k}} \nonumber}\\
&=[\alpha]\circ{\rm{Eval}}_{{{X_k\rightarrow \alpha}}} \circ S^{(l)}_{1,x_1,X_{k}}\circ_{u=1}^{{\lambda}}{\rm{Prim}}_{y_u,X_{k}}\nonumber\\
&=\sum_{h=0}^lb_h [\alpha]\circ{\rm{Eval}}_{{{X_k\rightarrow \alpha}}}  {\Big(\circ_{i=1}^hS_{1,{{\eta}}_i,X_{k}} \circ_{u=1}^{{\lambda}}{\rm{Prim}}_{y_u,X_{k}}\Big)\nonumber}\\
&=\sum_{h=0}^{{{\lambda}}-1}b_h [\alpha]\circ{\rm{Eval}}_{{{X_k\rightarrow \alpha}}}  {\Big(\circ_{u=1}^{{{\lambda}}-h}{\rm{Prim}}_{y_u,X_{k}}\Big)}+\sum_{h={{\lambda}}}^{l}b_h [\alpha]\circ{\rm{Eval}}_{{{X_k\rightarrow\alpha}}}  \circ S^{(l-h)}_{1,y_1,X_{k}}\nonumber\\
&=\sum_{h=0}^{{{\lambda}}-1}b_h \theta_{\alpha,{{\lambda}}-h,{{\lambda}}-h,X_{k}}+\sum_{h={{\lambda}}}^{l}b_h [\alpha]\circ{\rm{Eval}}_{{{X_k\rightarrow\alpha}}} \circ S^{(l-h)}_{1,y_1,X_{k}}\enspace. \label{tenkai zg}
\end{align}
We remark that, in the above equalities, we use $S_{1,y,X_{k}}\circ {\rm{Prim}}_{y,X_{k}}={\id}_{K[X_{k}]}$ for any {{$y\in K$}} which is not negative integer.
We apply Lemma~\ref{racinesymetrique2}, for
\begin{align*}  
\Delta_{\alpha}\circ \theta_{\alpha,{{\lambda}},{{\lambda}}-l,X_{k}}\circ \Theta_{\beta,\boldsymbol{l},\boldsymbol{l}}(P)&=\theta_{\alpha,{{\lambda}},{{\lambda}}-l,X_{k}}\circ \Theta_{\alpha,\boldsymbol{l},\boldsymbol{l}}\circ\Delta_{\alpha}(P)\enspace,
\end{align*} 
using $(\ref{tenkai zg})$, we obtain the assertion. {{\qed}}}}

{Let $l,k\in \Z$ with ${{0\leq  l \leq k}}$. 
We define a set of differential operators 
$$\mathcal{X}_{l,k}=\{V=\partial_1\circ \ldots \circ \partial_l\mid \partial_i\in \{1/\alpha, \tfrac{\partial}{\partial\alpha}-{{\rho}} y_1/\alpha\}, \ \# \{1\le i \le l , \partial_i=1/\alpha\}=k\}\enspace.$$

Let $\boldsymbol{{{\lambda}}}=({{\lambda}}_{i})_{1\le i \le {{\rho}}},\boldsymbol{{{\nu}}}=({{\nu}}_{i})_{1\le i \le {{\rho}}}$ be two integral vectors in $\zu^{{{\rho}}}$. By Lemma $\ref{derivative Theta}$, one gets that 
$$\frac{\partial^l}{\partial\alpha^l}(\Theta_{\alpha,\boldsymbol{{{\lambda}}},\boldsymbol{{{\nu}}}}({{\hat{P}}}))=\sum_{\mathbf{I}, \vert \mathbf{I}\vert\leq l}\sum_{V\in \mathcal{X}_{l,\vert \mathbf{I}\vert}}\Theta_{\alpha,\boldsymbol{{{\lambda}}},\boldsymbol{{{\nu}}}-\mathbf{I}}(V({{\hat{P}}}))\enspace.$$
By the Leibniz formula, $V({{\hat{P}}})$ is a linear combination (over $K[1/\alpha]$) of the derivatives 
$\frac{\partial^u}{\partial\alpha^u}({{\hat{P}}})$ for $0\leq {u}\leq l-\vert \mathbf{I}\vert$.
Since ${{\hat{P}}}=fg$, it is a linear combination of $g\frac{\partial^u}{\partial\alpha^u}(f)$, for $0\leq u \leq l-\vert\mathbf{I}\vert$. }
 
\

{{Relying on Lemma $\ref{racinesymetrique2}$ and Lemma $\ref{racinesymatrique3}$, together with the same argument of {\rm{\cite[Lemma~$9$]{DHK3}}}, we obtain}}
\begin{lemma} \label{sufficient condition 0}
$(${{\it confer}} {\rm{\cite[Lemma~$9$]{DHK3}}}$)$
Let $\boldsymbol{l}:=(1,2,\ldots,{{\rho}})\in \Z^{{{\rho}}}$ and $\bold{I}=(a_{1},\ldots,a_{{{\rho}}})\in \Z^{{{\rho}}}$ {{with $a_i\ge0$}}. 
Let $0\leq l$ be an integer with $\vert \bold{I}\vert \leq l$.
Assume further either of these two to be true~$:$
\begin{itemize} \label{condicombi}
\item[$({{a}})$] The $2{{\rho}}$ dimensional vector $(\boldsymbol{l},\boldsymbol{l}-\bold{I})$ has { two coordinates in common (in other words, either $i-a_i=j-a_j$ for some pair $i\neq j$ or for some $i$,  $1\leq i-a_i\leq \rho$, that is one of the coordinates of $\boldsymbol{l}$ is equal to one of the coordinates of $\boldsymbol{l}-\boldsymbol{I}$). }
\item[$({{b}})$] There exists an index $1\leq s \leq {{\rho}}$ such that $0\leq a_{s}-s<2{{\rho}} n-l+\vert\mathbf{I}\vert$.
\end{itemize}
Then, $\Delta_{\alpha}\circ \Theta_{\beta,\boldsymbol{l},\boldsymbol{l}}\circ\Theta_{\alpha,\boldsymbol{l},\boldsymbol{l}-\bold{I}}(g\frac{\partial^u f}{\partial\alpha^u})=0$ for all $0\leq u \leq l-\vert\mathbf{I}\vert$.
{Moreover, the smallest integer $l$
for which there exists $\mathbf{I}:=(a_{1},\ldots,a_{{{\rho}}})\in \Z^{{{\rho}}}$ satisfying  $a_i\ge0$} and $\vert\mathbf{I}\vert\leq l$  {where neither (a) nor (b) holds, is indeed $(2n+1){{\rho}}^2$. }
\end{lemma}

\smallskip

\smallskip

To completes the proof of Proposition $\ref{decompose Cnum}$, it remain to show $(iii)$.
{Since} Lemma $\ref{sufficient condition 0}$ ensures that  
$$\left.\frac{\partial^l}{\partial\alpha^l}(C_{u,m})\right|_{\alpha=\beta}=0 \hspace{15pt} \text{for all} \ 0\leq l\leq (2n+1){{\rho}}^2-1\enspace {,}$$
{it follows} that $C_{u,m}$ is {divisible} by $(\alpha-\beta)^{(2n+1){{\rho}}^2}$, which is $(iii)$.
Combining this result and Lemma $\ref{homo}$, we obtain Proposition $\ref{decompose Cnum}$.

\subsubsection{Fourth step}
We study the non-vanishing of the constant ${{c_{u,m}=}}c_{n,u,m}$ in {{Proposition}} $\ref{decompose Cnum}$. {At this stage we shall need to consider specializing the variables $t_{i,s_j}$ to $1$ instead of $\alpha_i$.  Since we are going to proceed by induction  on $m$, only one set of variables $t_{m,s_j}$ will be considered at a time, the others acting essentially as scalars. We shall thus denote by $\varphi_{1_i,x_j,s_j}$ the map taking $t_{i,s_j}^k\longmapsto {{1/(k+x_j+1)^{s_j}}}$ (compare with~$(\ref{defvarphi})$). Similarly, $\psi_{i,\boldsymbol{x}}={{\circ}}_{j=1}^d{{\circ}}_{s_{j}=1}^{r_{j}}\varphi_{1_i,x_{j},s_{j}}$, so that $\psi_{\boldsymbol{1},\boldsymbol{x}}={{\circ}}_{i=1}^m\psi_{i,\boldsymbol{x}}$. 
We shall write $\Xi={{\circ}}_{i=1}^{m-1}\psi_{i,\boldsymbol{x}}$, so that $\psi_{\boldsymbol{1},\boldsymbol{x}}=\Xi\circ\psi_{m,\boldsymbol{x}}$ (compare with~$(\ref{defpsi})$).}
\begin{lemma} \label{induction}
Set 
\begin{align*}
\displaystyle B_m(\boldsymbol{t}_m)=\displaystyle B(\boldsymbol{t})={\prod_{j=1}^d\prod_{s_j=1}^{r_j}} \left[t^{u}_{m,s_j} \cdot (t_{m,s_j}-{{1}})^{{{\rho}} n}\right] \prod_{(j_1,s_{j_1})<(j_2,s_{j_2})}(t_{m,s_{j_2}}-t_{m,s_{j_1}})\enspace.
\end{align*} 
Then, we have
$$c_{n,u,m}=(-1)^{{{\rho}}^2n(m-1)}c_{n,u+{{\rho}}(n+1),m-1}\cdot\psi_{m,\boldsymbol{x}}\left(B(\boldsymbol{t})\right)\enspace.$$
\end{lemma}

\proofname

Recall that by $(\ref{equality 1})$,
{\begin{align*} \label{deco Cnum2} 
D_{n,u,m}:&=\dfrac{C_{n,u,m}}{\prod_{i=1}^m\alpha^{{{\rho}}(u+1)+{{\rho}}^2n+\left(\sum_{1\le j_1<j_2\le d}r_{j_1}r_{j_2}+\sum_{j=1}^d{r_j\choose 2}\right)}_i}\\
&=c_{n,u,m}\prod_{1\le i_1<i_2\le m}(\alpha_{i_2}-\alpha_{i_1})^{(2n+1){{\rho}}^2}= \psi_{\boldsymbol{1},\boldsymbol{x}}\left(Q_{n,u,m}\right)\enspace. 
\end{align*}}
We are going to evaluate $D_{n,u,m}$ at $\alpha_m=0$  and thus separate the variables in $Q_{n,u,m}$ {(recall that $Q_{n,u,m}$ is defined in equation~$(\ref{defqu})$)} first. By definition, one has
\begin{equation*}
Q_{n,u,m}(\boldsymbol{t})=Q_{n,u,m-1}(\boldsymbol{t})\cdot B(\boldsymbol{t})C(\boldsymbol{t})\enspace,
\end{equation*} 
where 
\begin{align*}
C(\boldsymbol{t})&=\prod_{j=1}^{d}\prod_{s_j=1}^{r_j}\prod_{k=1}^{m-1}(\alpha_mt_{m,s_j}-\alpha_k)^{{{\rho}} n}\cdot\prod_{i=1}^{m-1} \prod_{j=1}^d\prod_{s_j=1}^{r_j} (\alpha_it_{i,s_j}-\alpha_m)^{{{\rho}} n}\cdot \prod_{i=1}^{m-1}\prod_{\substack{1\le j_1,j_2\le d \\ 1\le s_{j_1}\le r_{j_1} \\ 1\le s_{j_2}\le r_{j_2}}}\kern-8pt(\alpha_{m}t_{m,s_{j_2}}-\alpha_{i}t_{i,s_{j_1}})\enspace.
\end{align*}
Note that $Q_{n,u,m-1},B$ do not depend on $\alpha_m$, and {$\psi_{\boldsymbol{1},\boldsymbol{x}}$} treats $\alpha_m$ as a scalar.
Hence,
\begin{align} \label{compare 1}
\displaystyle  \ev_{\alpha_m\rightarrow 0}\left({D_{n,u,m}}\right) &= \displaystyle c_{n,u,m}\prod_{i=1}^{m-1}(-\alpha_i)^{(2n+1){{\rho}}^2}\prod_{1\le i<j\le m-1}(\alpha_{j}-\alpha_{i})^{(2n+1){{\rho}}^2}\enspace\\
&=\psi_{\boldsymbol{1}}\left(Q_{n,u,m-1}(\boldsymbol{t})B(\boldsymbol{t})\ev_{\alpha_m\rightarrow0}\left(C(\boldsymbol{t})\right)\right)\enspace. \nonumber
\end{align}
But { 
\begin{align*}
\ev_{\alpha_m\rightarrow0}\left(C(\boldsymbol{t})\right)&=\prod_{i=1}^{m-1}(-\alpha_i)^{{{\rho}}^2n}\prod_{i=1}^{m-1}\prod_{j=1}^d\prod_{s_j=1}^{r_j}(\alpha_it_{i,s_j})^{{{\rho}} n}
\prod_{i=1}^{m-1}\prod_{j=1}^{d}\prod_{s_j=1}^{r_j}(-\alpha_it_{i,s_j})^{{{\rho}}}\\
&=(-1)^{{{\rho}}^2(m-1)(n+1)}\prod_{i=1}^{m-1}\alpha_i^{(2n+1){{\rho}}^2}\prod_{i=1}^{m-1}\prod_{j-1}^d\prod_{s_j=1}^{r_j}t_{i,s_j}^{{{\rho}}(n+1)}\enspace.
\end{align*}}
We now note  that {$\psi_{m,\boldsymbol{x}}$} treats the variables {$t_{i,s_j}, 1\leq i\leq m-1$} as scalars and {$\Xi$} treats variables {$t_{m,s_j}$} as scalars and remark 
$$ 
Q_{n,u,m-1}(\boldsymbol{t})\ev_{\alpha_m\rightarrow0}(C(\boldsymbol{t}))=(-1)^{{{\rho}}^2(m-1)(n+1)}\prod_{i=1}^{m-1}\alpha_i^{(2n+1){{\rho}}^2}Q_{n,u+{{\rho}}(n+1),m-1}(\boldsymbol{t})\enspace.
$$
Thus
$$\begin{array}{lll}\displaystyle
\psi_{\boldsymbol{1},\boldsymbol{x}}\left(Q_{n,u,m-1}(\boldsymbol{t})B(\boldsymbol{t})\ev_{\alpha_m\rightarrow0}(C(\boldsymbol{t}))\right)\\ & 
\kern-90pt=&\rule{0mm}{7mm} \displaystyle(-1)^{{{\rho}}^2(m-1)(n+1)}\prod_{i=1}^{m-1}\alpha_i^{(2n+1){{\rho}}^2}\Xi(Q_{n,u+{{\rho}}(n+1),m-1}(\boldsymbol{t}))\psi_{m,\boldsymbol{x}}B(\boldsymbol{t}))\enspace.
\end{array}$$
Using the relation~$(\ref{compare 1})$, taking into account $D_{n,u+{{\rho}}(n+1),m-1}=\psi_{\boldsymbol{1},\boldsymbol{x}}(Q_{n,u+{{\rho}}(n+1),m-1}(\boldsymbol{t}))$ and simplifying,
$$c_{n,u,m}=
(-1)^{{{\rho}}^2n(m-1)}c_{n,u+{{\rho}}(n+1),m-1}\cdot \psi_{m,\boldsymbol{x}}(B(\boldsymbol{t}))\enspace.$$ 
This completes the proof of the lemma. \qed

{{By Lemma $\ref{induction}$, we have 
\begin{align} \label{A}
c_{n,u,m}&=(-1)^{{{\rho}}^2nm(m-1)/2}c_{n,u+(m-1){{\rho}}(n+1),1}\prod_{i=2}^m\psi_{i,\boldsymbol{x}}(B_i(\boldsymbol{t}_i))\enspace.
\end{align}
By the definition of $c_{n,u+(m-1){{\rho}}(n+1),1}$ (see Proposition $\ref{decompose Cnum}$), we have
\begin{align} \label{B}
c_{n,u+(m-1){{\rho}}(n+1),1}&={{\circ}}_{j=1}^d{{\circ}}_{s_j=1}^{r_j} \varphi_{\alpha_1,x_j,s_j}(B_{1}(\alpha_1,\boldsymbol{t}_1))\enspace,
\end{align}
where 
\begin{align*}
B_{1}(\alpha_1,\boldsymbol{t}_1)=
\prod_{j=1}^d \prod_{s_j=1}^{r_j} \left[t^{u+(m-1){{\rho}}(n+1)}_{1,s_j}(t_{1,s_j}-\alpha_1)^{{{\rho}} n}\right]\cdot \prod_{(j_1,s_{j_1})<(j_2,s_{j_2})}(t_{1,s_{j_2}}-t_{1,s_{j_1}})\enspace.
\end{align*}}}

\subsubsection{Last step}
{{Let $u$ {{be a non-negative integer}}. By $(\ref{A})$ and $(\ref{B})$, to complete the proof of non-vanishing of $c_{n,u,m}$, it is enough to prove
\begin{align} \label{C}
{{\circ}}_{j=1}^d{{\circ}}_{s_j=1}^{r_j} \varphi_{1,x_j,s_j}(B(\boldsymbol{t}))\neq 0\enspace,
\end{align}
where  
\begin{align*}
B(\boldsymbol{t})=
\prod_{j=1}^d \prod_{s_j=1}^{r_j} \left[t^{u}_{j,s_j}(t_{j,s_j}-1)^{{{\rho}} n}\right]\cdot \prod_{(j_1,s_{j_1})<(j_2,s_{j_2})}(t_{j_2,s_{j_2}}-t_{j_1,s_{j_1}})\enspace.
\end{align*}
To prove $(\ref{C})$, we {introduce the column vectors} 
$$M_h:=({{(s_j-1)!}}\varphi_{1, x_j,s_j}(t^{{{u}}+h{{-1}}}(t-1)^{{{\rho}} n}))_{1\le j \le d, 1\le s_j \le r_j}=\left(\int^1_0 t^{x_j+{{u}}+h{{-1}}}(t-1)^{{{\rho}} n}){\rm{log}}^{s_j-1}(1/t)\, dt\right)_{1\le j \le d, 1\le s_j \le r_j}\enspace,$$}} {{for $1\le h \le \rho$}} and let $M$ be the matrix whose columns are $M_h$.
Then we have $${{\circ}}_{j=1}^d{{\circ}}_{s_j=1}^{r_j} \varphi_{1,x_j,s_j}(B(\boldsymbol{t}))= {{\prod_{j=1}^d \prod_{s_j=1}^{r_j} \dfrac{1}{(s_j-1)!}}} \cdot {\rm{det}}(M)\enspace.$$
We now prove ${\rm{det}}(M)\neq 0$. {{By changing $x_j$ to $x_j+u$ for $1\le j \le d$, we may assume $u=0$.}}

We are going to prove the following key proposition.
\begin{proposition}
\label{integraleinf}
Suppose $x_i\neq x_j$ for all $i\neq j  \,(1\leq i\leq d)$. Then we have $\det\,(M)\neq 0$. 
\end{proposition}
We now study the matrix ${{M}}$. For this purpose, we denote by $C(x)$  {the row vector}:
{{$$C(x)=\left(\int_0^1t^{x+h-1}(t-1)^{{{\rho}} n}dt\right)_{1\leq h\leq {{\rho}}}\enspace.$$}

{{We have the following lemma.
\begin{lemma}\label{ecritureencolonneunevar}
$$M=\left(\frac{(-1)^{s_j-1}d^{s_j-1}C(x_j)}{dx_j^{s_j-1}}\right)_{1\leq j\leq d, 1\leq s_j\leq r_j}\enspace.$$
\end{lemma}}

\proofname

By simply noting that the integration with respect to $t$ and the differentiation with respect to $x$ commutes, one gets~:
$$\frac{d^k C(x)}{dx^k}=\left(\int_0^1\frac{d^k\left[t^{x+h-1}(t-1)^{{{\rho}} n}\right]}{d x^k}dt\right)=\left(\int_0^1t^{x+h-1}(t-1)^{{{\rho}} n}\log(t)^kdt\right)\enspace.$$
The lemma follows by adjusting the sign. \qed

\

\begin{lemma}\label{calcintfacile}
Let $y>0$ be a real number, and $m\geq 0$ be an integer. We have
$$F(y, m):=\int_0^1t^y(t-1)^m dt=\frac{(-1)^mm!}{\prod_{j=0}^{m}(y+j+1)}\enspace.$$
\end{lemma}

\proofname

{{The proof is simply done by induction on $m$ ({\it{confer}} Euler's Beta function in \cite{Askey} \cite{WW}), since
$$F(y, m)=\left[\frac{t^{y+1}}{y+1}\cdot (t-1)^{m}\right]_0^1-\frac{m}{y+1}\int_0^1t^{y+1}(t-1)^{m-1}=-\frac{m}{y+1}\cdot F(y+1,m-1)\enspace.$$ \qed}}

\begin{lemma}\label{Q}
Set ${{\tilde{Q}}}(x)=\prod_{{{j}}=1}^{{{\rho}}+{{\rho}} n}(x+{{j}})$ and, for $1\leq h\leq {{\rho}}$, 
${{\tilde{P}}}_h(x)=\prod_{j=1}^{h-1}(x+j)\prod_{j'={{h}}+1}^{{{\rho}}}(x+j'+{{\rho}} n)$
$($although these functions actually depend on ${{\rho}}, n$, we omit the subscripts, except in the proof below as they will be fixed and play no role$)$.
Then we have
$${}^t\!\kern1ptC(x)=(-1)^{{{\rho}} n}({{\rho}} n)!\cdot \left(\frac{{{\tilde{P}}}_h(x)}{{{\tilde{Q}}}(x)}\right)_{1\le h \le {{\rho}}}\enspace.\label{reformatfracratio}$$
\end{lemma}

\proofname

The lemma follows from Lemma~$\ref{calcintfacile}$, applied for $m={{\rho}} n$ and {$y=x+h-1$, $1\leq h\leq {{\rho}}$}.  \qed

We now set $\boldsymbol{P}(x)=({{\tilde{P}}}_h(x))_{1\leq h\leq {{\rho}}}$ and define the matrix $N$
$$N=\left(\frac{d^{s_j-1}\boldsymbol{P}(x_j)}{dx_j^{s_j-1}}\right)_{1\leq j\leq d, 1\leq s_j\leq r_j}\enspace.$$
We have
\begin{lemma}
$$\det(M)=\frac{(-1)^{{{\rho}}^2n-{{\rho}}/2+\left(\sum_{j=1}^dr_j^2\right)/2}({{\rho}} n)!^{{{\rho}}}}{\prod_{j=1}^d{{\tilde{Q}}}(x_j)^{{r_j}}}\det(N)\enspace.$$
\end{lemma}

\proofname

Taking into account Lemma~$\ref{ecritureencolonneunevar}$, we have  on the one hand~:
$$\det(M)=(-1)^{\sum_{j=1}^d\sum_{s_j=1}^{r_j}(s_j-1)}\det\left(\frac{d^{s_j-1}C(x_j)}{dx_j^{s_j-1}}\right)_{1\leq j\leq d, 1\leq s_j\leq r_j}\enspace.$$
On the other hand,
$$\det\left(\frac{d^{s_j-1}C(x_j)}{dx_j^{s_j-1}}\right)_{1\leq j\leq d, 1\leq s_j\leq r_j}=
\det\left(\frac{(-1)^{{{\rho}} n}({{\rho}} n)!d^{s_j-1}\frac{{{\tilde{P}}}_i(x_j)}{{{\tilde{Q}}}(x_j)}}{dx_j^{s_j-1}}\right)_{1\leq j\leq d, 1\leq s_j\leq r_j}\enspace.$$
Using the multi-linearity of determinant again,
$$\det\left(\frac{d^{s_j-1}C(x_j)}{dx_j^{s_j-1}}\right)_{1\leq j\leq d, 1\leq s_j\leq r_j}=(-1)^{{{\rho}}^2n}({{\rho}} n)!^{{{\rho}}}\det\left(\frac{d^{s_j-1}\frac{{{\tilde{P}}}_i(x_j)}{{{\tilde{Q}}}(x_j)}}{dx_j^{s_j-1}}\right)_{1\leq j\leq d, 1\leq s_j\leq r_j}\enspace.$$

However, since we have $$\frac{d\left(\frac{{{\tilde{P}}}_i(x_j)}{{{\tilde{Q}}}(x_j)}\right)}{dx_j}=\frac{\tfrac{d}{dx_j} {{\tilde{P}}}_i(x_j)}{{{\tilde{Q}}}(x_j)}-\frac{{{\tilde{P}}}_i(x_j)\tfrac{d}{dx_j} {{\tilde{Q}}}(x_j)}{{{\tilde{Q}}}^2(x_j)}\enspace,$$
the Leibniz formula and the multi-linearity of the determinant yields~:

$${\det(M)=\frac{(-1)^{\sum_{j=1}^d\sum_{s_j=1}^{r_j}(s_j-1)}\cdot (-1)^{{{\rho}}^2n}({{\rho}} n)!^{{{\rho}}}}{\prod_{j=1}^d{{\tilde{Q}}}(x_j)^{{r_j}}}\det(N)\enspace.}$$
The lemma follows up to the computation of the normalizing combinatorial factor. \\
To ease reading, we provide a couple of details~:
$$\sum_{j=1}^d\sum_{s_j=1}^{r_j}(s_j-1)=\sum_{j=1}^d\frac{r_j^2-r_j}{2}=\frac{\sum_{j=1}^dr_j^2-{{\rho}}}{2}\enspace,$$
and thus
$$(-1)^{\sum_{j=1}^d\sum_{s_j=1}^{r_j}(s_j-1)}\cdot (-1)^{{{\rho}}^2n}({{\rho}} n)!^{{{\rho}}}=(-1)^{{{\rho}}^2n-{{\rho}}/2+(\sum_{j=1}^dr_j^2)/2}({{\rho}} n)!^{{{\rho}}}\enspace. \ \ \qed$$ 

\bigskip

We now compute $\det(N)$. For this, we introduce a few notations.  We first define, for a set of variables $w_1, \ldots, w_\rho$,
$$V=V(w_1,\ldots, w_{{{\rho}}})=\det\left(w_j^i\right)_{\kern3pt0\leq i\leq {{\rho}}-1, 1\leq j\leq {{\rho}}}=\prod_{1\leq i<j\leq {{\rho}}}(w_j-w_i)\enspace,$$ the usual Vandermonde determinant and the degenerate version corresponding to Hermite interpolation problem~:
$$H=H(x_1,\ldots,x_d)=\left(\frac{d^{s_j-1}x_j^i}{dx_j^{s_j-1}}\right)_{0\leq i\leq {{\rho}}-1,1\leq j\leq d, 1\leq s_j\leq r_j}\enspace.$$
The evaluation of the Hermite determinant is indeed classical, found in T.~Muir's encyclopedia, \cite{muir}. It is attributed to L. Schendel, \cite{schen} and
recent appearances in the literature can be found for instance in \cite{vander} by A. van der Poorten and by Loring W. Tu  \cite{tu}.
We rely on this evaluation in the following lemma.
\begin{lemma} \label{H}
 One has
$$\det(H)=(-1)^{\sum_{j=1}^dr_j(r_j-1)/2}\cdot \left(\prod_{i=1}^{d}\prod_{s_i=0}^{r_i-1}s_i!\right)\cdot\prod_{1\leq i<j\leq d}(x_j-x_i)^{r_ir_j}\enspace.$$
\end{lemma}
We now compute special values of the particular interpolation determinant we are interested in.
\begin{lemma} \label{N tilde}
 Define the non-degenerate analog of $N$ as
$${{\tilde{N}(w_1,\ldots,w_{\rho})}}=\left(\boldsymbol{P}(w_j)\right)_{1\leq j\leq {{\rho}}}\enspace.$$
Then,
$$\det(\tilde{N}(-1,\ldots,-{{\rho}}))=(-1)^{{{\rho}}({{\rho}}-1)/2}\prod_{i=1}^{{{\rho}}-1}(i!)^2\prod_{i=1}^{{{\rho}}-1}{i+{{\rho}} n\choose {{\rho}} n}\enspace.$$
\end{lemma}

\proofname

Recall that ${{\tilde{P}}}_i(x)=\prod_{j=1}^{i-1}(x+j)\prod_{j={{i}}+1}^{{{\rho}}}(x+j+{{\rho}} n)$, so, if $-{{\rho}} \leq -i\leq -1$,
one has $$\boldsymbol{P}(-i)={}^t{{\begin{pmatrix} \star ,\ldots, \star, &  {{\tilde{P}}}_i(-i), & 0, & \ldots, & 0 \end{pmatrix}}}\enspace.$$
It follows that 
$$\det\left(\tilde{N}(-1,\ldots,-{{\rho}})\right)=\prod_{i=1}^{{{\rho}}}{{\tilde{P}}}_i(-i)=(-1)^{{{\rho}}({{\rho}}-1)/2}\prod_{i=1}^{{{\rho}}}(i-1)!\prod_{i=1}^{{{\rho}}-1}i!\prod_{i=1}^{{{\rho}}-1}{i+{{\rho}} n\choose {{\rho}} n}\enspace.$$ \qed

We are now in a position to compare determinants.
\begin{lemma} \label{normaconstante} 
We have the following equality between determinants. 
Let $\boldsymbol{R}$ be the matrix representing the polynomials {$({{\tilde{P}}}_i(x))_{1\le i \le {{\rho}}}$ in the standard basis ${}^t(1,\ldots,x^{{{\rho}}-1})$}. This is a square ${{\rho}}\times {{\rho}}$ matrix.
Then, 
$$\det\left(\tilde{N}(w_1,\ldots,w_{{{\rho}}})\right)= \det(\boldsymbol{R})\cdot V\enspace,$$
and
$$\det\left({N}(x_1,\ldots,x_d)\right)= \det(\boldsymbol{R})\cdot \det(H)\enspace.$$
\end{lemma}

\proofname

{{The matrix $N$ is obtained from the matrix $H$ by multiplication by $\boldsymbol{R}$ and similarly, $\tilde{N}$ is obtained by multiplying the Vandermonde matrix by $\boldsymbol{R}$. \qed}}

\begin{lemma} \label{last} 
One has
$$\det\left({N}(x_1,\ldots,x_d)\right)=(-1)^{\sum_{j=1}^dr_j(r_j-1)/2}\cdot \frac{\prod_{i=1}^{{{\rho}}-1}(i!)^2\prod_{i=1}^{{{\rho}}-1}{i+{{\rho}} n\choose {{\rho}} n}}{\prod_{1\leq i<j\leq {{\rho}}}(j-i)}\cdot \prod_{i=1}^{d}\prod_{s_j=0}^{r_j-1}s_j!\cdot \prod_{1\leq i<j\leq d}(x_j-x_i)^{r_ir_j}\enspace.$$
\end{lemma}

\proofname

It is enough to use Lemma~$\ref{normaconstante}$ with the Vandermonde equality at $(-1,\ldots,-{{\rho}})$ to compute $\det(\boldsymbol{R})$ and then plug in the constant in the equality linking $\det(N)$ with the determinant of the Hermite matrix of Lemma~$\ref{normaconstante}$. One finds
$$\det(\boldsymbol{R})=\frac{\det(\tilde{N}(-1,\ldots,-{{\rho}}))}{V(-1,\ldots,-{{\rho}})}=\frac{(-1)^{{{\rho}}({{\rho}}-1)/2}\prod_{i=1}^{{{\rho}}-1}(i!)^2\prod_{i=1}^{{{\rho}}-1}{i+{{\rho}} n\choose {{\rho}} n}}{\prod_{1\leq i<j\leq {{\rho}}}(i-j)}=\frac{\prod_{i=1}^{{{\rho}}-1}(i!)^2\prod_{i=1}^{{{\rho}}-1}{i+{{\rho}} n\choose {{\rho}} n}}{\prod_{1\leq i<j\leq {{\rho}}}(j-i)}\enspace.$$
Thus, by above equality and Lemma $\ref{H}$, we obtain the desire equality.  \qed

\section{Proof of Theorems}
\label{estim}
\subsection{{{Analytic estimates in Pad\'e approximation}}}
In this section, we use the following notations. Let $K$ be a number field and $v$ a place of $K$. {{Denote by}} $K_v$ the completion of $K$ at $v$ and $\vert\cdot\vert_v$ by the absolute value corresponding to $v$. Throughout the section, the small $o$-symbol  $o(1)$ refers when $n$ tends to infinity. 
Put $\varepsilon_v=1$ if $ v|\infty$ and $0$ otherwise.

Let $I$ be a non-empty finite set of indices, $A=K[\alpha_i]_{i\in I}[z,t]$ be a polynomial ring in indeterminate $\alpha_i,z,t$. 
We set $\Vert P\Vert_v=\max\{\vert a\vert_v\}$ where $a$ runs in the coefficients of $P$. Thus the ring $A$ is endowed with a structure of normed vector space. If $\phi$ is an endomorphism of $A$, we denote by $\Vert \phi\Vert_{v}$ the endomorphism norm defined in a standard way $\Vert \phi\Vert_v=\inf\{M\in\ru, \forall\kern3pt x\in A, \kern3pt\Vert\phi(x)\Vert_v\leq M\Vert x\Vert_v\}=\sup\left\{\frac{\Vert\phi(x)\Vert_v}{\Vert x\Vert_v}, 0\neq x\in A\right\}$. This norm is well defined provided $\phi$ is continuous. Unfortunately, we will have to deal  also with non-continuous morphisms. In such a situation, we restrict the source space to some appropriate sub-vector space $E$ of $A$ and talk of $\Vert\phi\Vert_v$ with $\phi$ seen  as $\left.\phi\right|_{E}: E\longrightarrow A$ on which $\phi$ is continuous. In case of perceived ambiguity, it will be denoted by $\Vert\phi\Vert_{E,v}$.
The degree of an element of $A$ is as usual the total degree {{with respect to all the $\alpha_i, t$ and $z$}}. 
\begin{lemma}[{{\it confer}} \cite{DHK3}, Lemma~13] \label{(vi)}
Let $E_N$ be the subspace of $A$ consisting of polynomials of degree at most $N$ in $t$. Let $x=a/b\in \Q\cap [0,1)$ with $a,b$ are coprime integers.
Then $\inte_x:E_N\longrightarrow A$ as defined in Notation~$\ref{notationderiprim}$ satisfies 
$$\Vert \inte_{{x}}\Vert_{E_N,v}\leq \max\{\left\vert 1/(k+x+1)\right\vert_v,0\leq k\leq N\}\leq \vert d_{N+1}(a,b) \vert_v^{-{ (1-\varepsilon_v)}}\enspace,$$ with $d_N(a,b)={\rm{l.c.m.}}(a,a+b,\ldots,a+b{{N}})$ {{where}} {\rm{l.c.m.}} denotes the least common multiple.
\end{lemma}
From the preceding lemma, we deduce
the following estimates.
\begin{lemma}[{{\it confer}} \cite{DHK3}, Lemma~14] \label{majonorme} 
Let $x_j=a_j/b_j\in \Q\cap [0,1)$ for $1\le j \le d$ {{as in Theorem $\ref{Lerch}$}}. We denote by $w$ the place of $\Q$ such that $v\vert w$. Then, one has~$:$
\begin{itemize}
\item[$(i)$] The polynomial $P_l(z)=P_{n,l}(\boldsymbol{\alpha},\boldsymbol{x}\vert z)$ satisfies 
$$\log \Vert P_{l}(z)\Vert_v\leq 
\begin{cases}
\dfrac{n{{[K_v:\Q_w]}}}{{{[K:\Q]}}}\left({{\rho}} m\log(2)+{{\rho}}\left(\log({{\rho}} m+1)+{{\rho}} m\log\left(\dfrac{{{\rho}} m+1}{{{\rho}} m}\right)\right)+{o(1)}\right) & \  \text{if} \ v|\infty\\
\displaystyle{\sum_{j=1}^d} r_j \log\, |\mu_n(x_j)|^{-1}_v & \  \text{otherwise} \enspace.
\end{cases}
$$ Recall that $P_{l}(z)$ is of degree ${{\rho}} n$ in each variable $\alpha_i$, of degree ${{\rho}} mn+l$ in $z$ and constant in $t$.
\item[$(ii)$] 
The polynomial $P_{l,i,j,s_j}(z)=P_{n,l,i,j,s_j}(\boldsymbol{\alpha},\boldsymbol{x} \vert z)$ satisfies
\begin{align*}
\log \Vert P_{l,i,j,s_j}(z)\Vert_v&\leq  
\begin{cases}
\dfrac{n{{[K_v:\Q_w]}}}{{{[K:\Q]}}}\left({{\rho}}m\log(2)+\rho\left(\log({{\rho}} m+1)+{{\rho}} m\log\left(\dfrac{{{\rho}} m+1}{{{\rho}} m}\right)\right)+{o(1)}\right)& \  \text{if} \ v|\infty\\
\displaystyle{\sum_{j=1}^d} r_j \log\, |\mu_n(x_j)|^{-1}_v+{{\log \left \vert d_{{{\rho}} m(n+1)+1}(a_j,b_j)\right \vert_v^{-s_j}}} & \  \text{otherwise} \enspace.
\end{cases}
\end{align*}
Also, $P_{l,i,j,s_j}(z)$ is of degree $\leq {{\rho}} mn+l$ in $z$, of degree ${{\rho}} n$ in each of the variables $\alpha_j$ except for the index $i$ where it is of degree ${{\rho}} n+1$ $($recall that $\varphi_{\alpha_i,x_j,s_j}$ involves multiplication by 
$[\alpha_i]$$)$.
\item[$(iii)$] For any integer $k\geq 0$, the polynomial $\varphi_{\alpha_i,x_j,s_j}\circ[t^{k+n}](P_{l}(t))$ satisfies 
{{\small{\begin{align*}
\log \Vert \varphi_{\alpha_i,x_j,s_j}\circ[t^{k+n}](P_{l}(t)) \Vert_v \leq
\begin{cases}
\dfrac{n{{[K_v:\Q_w]}}}{{{[K:\Q]}}}\left({{\rho}} m\log(2)+{{\rho}} \left(\log({{\rho}} m+1)+{{\rho}} m\log\left(\dfrac{{{\rho}} m+1}{{{\rho}} m}\right)\right)+{o(1)}\right)& \  \text{if} \ v|\infty\\
\displaystyle{\sum_{j=1}^d} r_j \log\, |\mu_n(x_j)|^{-1}_v+{{\log \left \vert d_{{{\rho}} m(n+1)+1}(a_j,b_j)\right \vert_v^{-s_j}}} & \  \text{otherwise} \enspace.
\end{cases}
\end{align*}
}}}
By definition, it is a homogeneous polynomial in just the variables $\boldsymbol{\alpha}$ of degree $\leq {{\rho}} mn+l+k+n+1$.
\end{itemize}
\end{lemma}

\proofname

For a non-negative integer $N$, we denote $E_N$ the sub-vector space of $K[y_1,\ldots,y_m, z,t]$ consisting of polynomials of degree at most $N$ in the variables $y_i$ and define a morphism 
$$\Gamma: E_N\longrightarrow K[\alpha_1,\ldots,\alpha_m,z,t]; \ P(y_1,\ldots,y_m,z,t)\mapsto P(t-\alpha_1,\ldots,t-\alpha_m,z,t)\enspace.$$ 
Then, by {\rm{\cite[Lemma $5.3$ $(iv)$]{DHK3}}}, we have $$\Vert\Gamma \Vert_{E_N,v}\leq \left(2^{Nm}(N+1)^{m}\right)^{\varepsilon_v{ {{[K_v:\Q_w]}}/{{[K:\Q]}}}}\enspace.$$

\bigskip

Set $B_{n,l}(\boldsymbol{y},t)=t^l\prod_{i=1}^my_i^{{{\rho}} n}$, since $B_{n,l}$ is a monomial, its norm $\Vert B_{n,l}\Vert_v=1$. 
{By definition,  $P_{l}(z)=\ev_{{{t \rightarrow z}}}\circ_{j=1}^d S^{(r_j)}_{n,x_j}\circ\Gamma(B_{n,l})$}, and thus, by {{sub-multiplicativity}} of the endomorphism norm, 
$$\Vert P_{l}( z)\Vert_v\leq 2^{\varepsilon_v {{\rho}} mn  {{[K_v:\Q_w]}}/{{[K:\Q]}}} \prod_{j=1}^d \left\Vert S_{n,x_j} \right\Vert_v^{r_j}{{e^{\varepsilon_v n{o(1)}}}}\enspace,$$ (one can use  {\rm{\cite[Lemma $5.2$ $(iv)$]{DHK3}}} and 
$N={{\rho}} m(n+1)+l$ for property {\rm{\cite[Lemma $5.2$ (v)]{DHK3}}} using $0\leq l\leq {{\rho}} m$, and note that the original polynomial is a constant in $z$ so the evaluation map is an isometry). 

\vspace{0.5\baselineskip}

\bigskip

Now, by definition, $P_{l,i,j,s_j}(z)=\varphi_{\alpha_i,x_j,s_j}\circ \Theta(P_{l}(z))$ where
$$\Theta: K[\alpha_1,\ldots,\alpha_m,z,t] \longrightarrow K[\alpha_1,\ldots,\alpha_m,z,t]; \ P\mapsto \frac{P(\alpha_i,z)-P(\alpha_i,t)}{z-t}\enspace.$$
By definition, $\varphi_{\alpha_i,x_j,s_j}=[\alpha_i]\circ \ev_{{{t\rightarrow \alpha_i}}}\circ\inte^{(s_j)}_{x_j}$. 
Using again {\rm{\cite[Lemma $5.2$ $(iii),(vi)$]{DHK3}}} with $N={{\rho}} m(n+1)$ and {\rm{\cite[Lemma $5.2$ $(i),(ii)$]{DHK3}}}, with $N={{\rho}} m(n+1)$, and since ${{\rho}} mn=n\exp({o(1)})$, one gets $(ii)$.

Finally, {{we have}} $\varphi_{\alpha_i,x_j,s_j}\circ[t^{k+n}](P_{l}(t))=[\alpha_i]\circ\ev_{{{t\rightarrow \alpha_i}}}\circ \inte^{(s_j)}_{x_j}\circ
[t^{k+n}](P_{l}(t))$.
Again, using {\rm{\cite[Lemma $5.2$]{DHK3}}}, one gets $(iii)$. \qed

Recall that if $P$ is a homogeneous polynomial in some variables $y_i,i\in I$, for any point $\boldsymbol{\alpha}=(\alpha_i)_{i\in I}\in K^{\card(I)}$ where $I$ is any finite set, and $\Vert\cdot\Vert_v$ stands for the sup norm in $K_v^{\card(I)}$, with $$C_v(P)= (\deg(P)+1)^{\frac{\varepsilon_v{{[K_v:\Q_w]}}(\card(I))}{{{[K:\Q]}}}}\enspace,$$
one has
\begin{equation}\label{estimhomo}
\vert P(\boldsymbol{\alpha})\vert_v\leq C_v(P) \Vert P\Vert_v\cdot {\Vert {\boldsymbol{\alpha}\Vert_v^{\deg(P)}}}\enspace.
\end{equation}
So, the preceding lemma yields trivially estimates for the $v$-adic norm of the above given polynomials.
\begin{remark} It may be worthwhile to note that it is also possible to estimate with respect to ${\mathrm{h}}_v(\alpha_i-\beta)$ instead of ${\mathrm{h}}_v(\alpha_i)$ which can be useful in certain circumstances (saving of the archimedean error term ${{\rho}} mn\log(2)$, useful if the local heights of $\beta,\alpha_i$ are not too far apart).
\end{remark}
\begin{lemma} [{{\it confer}} \cite{DHK3}, Lemma~15] \label{upper jyouyonew}
Let $n$ be a positive integer, {{$x_j={{a_j/b_j}}\in \Q\cap [0,1)$}} and $\beta\in K$ with $\Vert \boldsymbol{\alpha}\Vert_v<\vert \beta\vert_v$.
{{We denote by $w$ the place of $\Q$  such that  $v\vert w$.}} 
Then we have for all $0\leq l\leq {{\rho}} m, {{1 \leq}} i\leq m, {{1\le j \le d}}$ and ${{1\le s_j\leq r_j}}$,
\begin{align*}
\log |R_{l,i,j,s_j}(\beta)|_v &\le {{{\rho}} m(n+1)}\log \Vert \boldsymbol{\alpha}\Vert_v
+(n+1)\log \left(\frac{\Vert\boldsymbol{\alpha}\Vert_v}{\vert \beta\vert_v}\right)+\log \left(\frac{\varepsilon_v\vert \beta\vert_v}{\vert \beta\vert_v-\Vert \boldsymbol{\alpha}\Vert_v}+(1-\varepsilon_v)\right)\\
&{{+}}\begin{cases}
\dfrac{n{{[K_v:\Q_w]}}}{{{[K:\Q]}}}\left({{\rho}} m\log(2)+{{\rho}}\left(\log({{\rho}} m+1)+{{\rho}} m\log\left(\dfrac{{{\rho}} m+1}{{{\rho}} m}\right)\right)+{o(1)}\right)& \  \text{if} \ v|\infty\\
\displaystyle{\sum_{j=1}^d} r_j \log\, |\mu_n(x_j)|^{-1}_v+{{\log \left \vert d_{{{\rho}} m(n+1)+1}(a_j,b_j)\right \vert_v^{-s_j}}} & \  \text{otherwise} \enspace.
\end{cases}
\end{align*}
\end{lemma} 

\proofname

By the definition of $P_{l}(z)$, as formal power series, we have
$$
R_{l,i,j,s_j}(z)=\sum_{k=0}^{\infty}\frac{\varphi_{\alpha_i,x_j,s_j}(t^{k+n}P_{l}(t))}{z^{k+n+1}} \enspace, 
$$
using the triangle inequality, the fact that $0\leq l\leq {{\rho}} m$ and Lemma~$\ref{majonorme}$, $(iii)$, together with inequality $(\ref{estimhomo})$
{\small{\begin{align*}
\vert R_{l,i,j,s_j}(\beta)\vert_v&\leq \Vert \boldsymbol{\alpha}\Vert_v^{{{\rho}} m(n+1)}\sum_{k=0}^{\infty}\left(\frac{\Vert\boldsymbol{\alpha}\Vert_v}{\vert \beta\vert_v}\right)^{n+1+k}
\cdot  2^{\varepsilon_v {{\rho}} mn {{[K_v:\Q_w]}}/{{[K:\Q]}}} \prod_{j=1}^d \left\Vert S_{n,x_j} \right\Vert_v^{r_j}\cdot \kern-2pt
\begin{cases}
e^{n{o(1)}} & \text{if} \ v|\infty \\
{{\left\vert d_{{{\rho}} m(n+1)+1}(a_j,b_j)\right\vert_v^{-s_j}}} & \text{otherwise}\enspace,
\end{cases}
\end{align*}}}
and the lemma follows using geometric series summation. \qed

\subsection{Proof of Theorem $\ref{Lerch}$}
To prove Theorem $\ref{Lerch}$, we show the following theorem. 
\begin{theorem} \label{Lerch 2}
We use the same notations as in Theorem $\ref{Lerch}$. {{Denote $x_j=a_j/b_j$ with $a_j,b_j\in\Z$ coprime.}} 
For a {{strictly positive integer}} $n$, we put $$D_n(m,\boldsymbol{x})=D_n={\rm{l.c.m.}}(a_j,a_j+b_j,\ldots,a_j+b_j({{\rho}} m(n+1)+1)\mid j=1,\ldots,d)\enspace.$$
For any place  $v\in \mathfrak{M}_K$, we define the constants
\begin{align*}
c(\boldsymbol{x},v)=\varepsilon_{v}{{\frac{{{[K_v:\Q_w]}}}{{{[K:\Q]}}}}}\left[\const\right]{{+}}(1-\varepsilon_{v})\sum_{j=1}^dr_j\log\, {|\mu(x_j)|_v^{{{-1}}}} \enspace,
\end{align*}
where $\varepsilon_v=1$ if $ v|\infty$ and $0$ otherwise.
We also define
\begin{align*}
\mathbb{A}(\boldsymbol{\alpha},\boldsymbol{x},\beta)&=\log\vert\beta\vert_{v_0}-({{\rho}} m+1)\log\Vert \boldsymbol{\alpha}\Vert_{v_0}-c(\boldsymbol{x},v_0)+(1-\varepsilon_{v_0})\cdot \max_{1\le j \le d}(r_j)\lim_{n\to \infty} \dfrac{\log\,|D_n|_{v_0}}{n} \enspace\kern-3pt,\\
U(\boldsymbol{\alpha},\boldsymbol{x},\beta)&={{\rho}} m{\mathrm{h}}_{v_0}(\boldsymbol{\alpha},\beta)+c(\boldsymbol{x},v_0)\enspace\kern-3pt,
\end{align*}
and finally 
\begin{align*}
V(\boldsymbol{\alpha},\boldsymbol{x},\beta)&={\rm{log}}\,|\beta|_{v_0}-{{\rho}} m{\mathrm{h}}(\boldsymbol{\alpha},\beta)-{{\rho}} m\log\,\|\boldsymbol{\alpha}\|_{v_0}+{{\rho}} m\log\,\|(\boldsymbol{\alpha},\beta)\|_{v_0}
-\sum_{j=1}^dr_j\log\,\mu(x_j)\\
&-\left[\const\right]-\max_j(r_j)\lim_{n\to \infty}\dfrac{\sum_{v\in \mathfrak{M}^f_K}\log\,|D_n|_v}{n}\enspace.
\end{align*}
Assume $V(\boldsymbol{\alpha},\boldsymbol{x},\beta)>0$. 
Then for any positive number $\varepsilon$ with $\varepsilon<V(\boldsymbol{\alpha},\boldsymbol{x},\beta)$, there exists an effectively computable positive number $H_0$ depending on $\varepsilon$ and the given data such that the following property holds.
For any $\boldsymbol{\lambda}:=(\lambda_0,\lambda_{i,j,s_j})_{\substack{1\le i \le m \\ 1\le j \le d, 1\le s_j \le r_j}} \in K^{{{\rho}} m+1} \setminus \{ \bold{0} \}$ satisfying $H_0\le {\mathrm{H}}(\boldsymbol{\lambda})$, then we have
\begin{align*}
\left|\lambda_0+\sum_{i=1}^m\sum_{j=1}^d\sum_{s_j=1}^{r_j}\lambda_{i,j,s_j}\Phi_{s_j}(x_j,\alpha_i/\beta)\right|_{v_0}
>C(\boldsymbol{\alpha},\boldsymbol{x},\beta,\varepsilon) {\mathrm{H}}_{v_0}(\boldsymbol{\lambda}){\mathrm{H}}(\boldsymbol{\lambda})^{-\mu(\boldsymbol{\alpha},\boldsymbol{x},\beta,\varepsilon)}\enspace,
\end{align*}
where  
\begin{align*}
&\mu(\boldsymbol{\alpha},\boldsymbol{x},\beta,\varepsilon):=
\dfrac{\mathbb{A}(\boldsymbol{\alpha},\boldsymbol{x},\beta)+U(\boldsymbol{\alpha},\boldsymbol{x},\beta)}{V(\boldsymbol{\alpha},\boldsymbol{x},\beta)-\varepsilon}\enspace, \\
&C(\boldsymbol{\alpha},\boldsymbol{x},\beta,\varepsilon):=\exp\left(-\left(\frac{\log(2)}{V(\boldsymbol{\alpha},\boldsymbol{x},\beta)-\varepsilon}+1\right)
\left(\mathbb{A}(\boldsymbol{\alpha},\boldsymbol{x},\beta)+U(\boldsymbol{\alpha},\boldsymbol{x},\beta)\right)\right)\enspace.
\end{align*}
\end{theorem}

\proofname

By Proposition $\ref{non zero det}$, the matrix ${{\mathrm{M}}}_n=\begin{pmatrix}
P_{l}(\beta)\\
P_{l,i,j,s_j}(\beta)
\end{pmatrix}$ 
with entries in $K$ is invertible. We apply a linear independence criterion in \cite[Proposition $5.6$]{DHK3}. 
We have by Lemma~$\ref{majonorme}$, $(i)$ together with inequality $(\ref{estimhomo})$
\begin{align*}
\log\Vert  P_{l}(\beta)\Vert_{v} & \leq \displaystyle \varepsilon_v\left({\frac{n{{\rho}} {{[K_v:\Q_w]}}}{{{[K:\Q]}}}}\left[\const+{o(1)}\rule{0mm}{4mm}\right]\right)\\
&+(1-\varepsilon_v)\sum_{j=1}^dr_j\log\left\vert\mu_n(x_j)\right\vert_v^{-1}+({{\rho}} mn+l){\mathrm{h}}_{v}(\boldsymbol{\alpha},\beta)\\ 
&\leq \displaystyle n\left({{\rho}} m{\mathrm{h}}_{v}(\boldsymbol{\alpha},\beta)+c(\boldsymbol{x},v)\right)+{o(1)}=n\displaystyle {{U}}(\boldsymbol{\alpha},\boldsymbol{x},\beta)+{o(1)}\enspace.
\end{align*}
Similarly, using this time Lemma~$\ref{majonorme}$ $(ii)$ {{and inequality $(\ref{estimhomo})$}},
\begin{align*}
\log\Vert  P_{l,i,j,s_j}(\beta)\Vert_{v} & \leq \displaystyle \varepsilon_v\left({{\frac{n{{\rho}} {{[K_v:\Q_w]}}}{{{[K:\Q]}}}}}\left[\const+{o(1)}\rule{0mm}{4mm}\right]\right)\\
&+(1-\varepsilon_v)\log\left\vert\mu_n(x)\right\vert_v^{-1}+({{\rho}} mn+l){\mathrm{h}}_{v}(\boldsymbol{\alpha},\beta)-(1-\varepsilon_v)s_j\log\vert D_n \vert_v\\ 
& \leq  \displaystyle n\left({{\rho}} m{\mathrm{h}}_{v}(\boldsymbol{\alpha},\beta)+c(\boldsymbol{x},v)\right)+f_v(n)\enspace,
\end{align*}
where
\begin{align*}
f_v:\N\rightarrow \R_{\ge0}; \ \ n\longmapsto  {{\rho}} m{\rm{h}}_v(\boldsymbol{\alpha},\beta)-(1-\varepsilon_v)\max_{1\le j\le d} \,(r_j) \cdot \log\,|D_n|_v\enspace.
\end{align*}
We define $$\displaystyle {{F_v}}(\boldsymbol{\alpha},\boldsymbol{x},\beta):\N\rightarrow \R_{\ge0}; \ n\mapsto n\left({{\rho}} m{\mathrm{h}}_{v}(\boldsymbol{\alpha},\beta)+c(\boldsymbol{x},v)\right)+f_v(n)\enspace.$$
By Lemma~$\ref{upper jyouyonew}$ ensures 
$$\log\, \vert  R_{l,i,j,s_j}(\beta)\vert_{v_0}\leq n\left(-\mathbb{A}(\boldsymbol{\alpha},\boldsymbol{x},\beta)+{o(1)}\right)\enspace,$$
by a linear independence criterion in \cite[Proposition $5.6$]{DHK3} for $\{\theta_{i,j,s_j}:=\Phi_{s_j}(x_j,\alpha_i/\beta)\}_{\substack{1\le i \le m\\ 1\le j \le d, 1\le s_j \le r_j}}$ and the above data, 
we obtain the assertions of Theorem $\ref{Lerch 2}$. 
Using the prime number theorem, we have 
$${\displaystyle{\lim_{n\to \infty}}}\dfrac{{\sum_{v\in \mathfrak{M}_K}}f_v(n)}{n}=\max_{j}(r_j)\lim_{n\to \infty}\dfrac{\sum_{v\in \mathfrak{M}^f_K}\log\,|D_n|_{v}}{n}\le \max_{j}(r_j)b{{\rho}} m\enspace,$$ 
and
$$
\sum_{v\in \mathfrak{M}_K}c(\boldsymbol{x},v)=\const+\sum_{j=1}^dr_j{\rm{log}}\, \mu(x_j)\enspace,
$$
we conclude
\begin{align*}
\mathbb{A}(\boldsymbol{\alpha},\boldsymbol{x},\beta)-\lim_{n\to \infty}\dfrac{1}{n}\sum_{v\neq v_0}{{F_v}}(\boldsymbol{\alpha},\boldsymbol{x},\beta)(n)&\ge  {\rm{log}}\, |\beta|_{v_0}-{{\rho}} m{\mathrm{h}}(\boldsymbol{\alpha},\beta)
+{{\rho}} m{\rm{log}}\,\|(\boldsymbol{\alpha},\beta)\|_{v_0}-({{\rho}} m+1)\log\,\|\boldsymbol{\alpha}\|_{v_0}\\
&-\sum_{v\in \mathfrak{M}_K}c(\boldsymbol{x},v)-\max_{j}(r_j)\lim_{n\to \infty}\dfrac{\sum_{v\in \mathfrak{M}^f_K}\log\,|D_n|_{v}}{n}
\ge V(\boldsymbol{\alpha},\boldsymbol{x},\beta)\enspace. \ \ \qed
\end{align*} 

{\subsection{Proof of Theorem~$\ref{periodiclerch}$}
We show that the periodic case (Theorem~$\ref{periodiclerch}$) is a corollary of Theorem~$\ref{Lerch}$.
\begin{lemma} 
In the situation of Theorem~$\ref{periodiclerch}$, for $x\in\Q$ which is not negative integer and  any integer $s\geq 1$, the $K$-subspace of $K_{v_0}$ generated by
$\{f_{b,w_i,x,s}(\beta)\mid 0 \le i \le m-1\}$ is the one generated by $\{\Phi_{s}(x,\alpha_i/\beta)\mid 0\le i \le m-1\}$.
\label{elementsimple}
\end{lemma}

\proofname

{{We denote the former $K$-subspace by $V_1$ and the latter by $V_2$.}}
Let $w(z)\in K[z]$ and assume ${\rm{deg}}\,w(z)\le q-1$. Then the rational function $w(z)/b(z)$ can be written as
$$\frac{w(z)}{{{b}}(z)}=\sum_{i=1}^m\frac{\gamma_i}{z-\alpha_i}\enspace,$$
for some $\gamma_i\in K$. 
Hence, 
$$\frac{w(z)}{{{b}}(z)}=\sum_{i=1}^m\gamma_i\sum_{k=0}^{\infty}\frac{\alpha_i^k}{z^{k+1}}\enspace.$$
By definition, 
$$f_{b,w,x,s}(z)=\sum_{i=1}^{m}\gamma_i\Phi_{s}(x,\alpha_i/z)\enspace.$$
{Above identity yields $V_1\subseteq V_2$.}
Conversely, since $w_0,\ldots,w_{m-1}$ are linearly independent, the rational functions ${{1/(z-\alpha_i)}}$ are linear combinations of  ${w_0(z)}/{b(z)}, \ldots,{w_{m-1}(z)}/{b(z)}$. 
By the same argument as above, we obtain $V_2\subseteq V_1$. \qed

Now, using Theorem~$\ref{Lerch}$ we know that provided $V(\boldsymbol{\alpha},\boldsymbol{x},\beta)>0$, the vector space  $$K\cdot 1 + \sum_{j=1}^{d}\sum_{s_j=1}^{r_j}\sum_{i=1}^mK\cdot \Phi_{s_j}(x_j,\alpha_i/\beta)$$ is of dimension ${{\rho}} m+1$ and is equal to $$K\cdot 1 + \sum_{j=1}^{d}\sum_{s_j=1}^{r_j}\sum_{i=1}^mK\cdot f_{b,w_i,x_j,s_j}(\beta)$$ by Lemma~$\ref{elementsimple}$. This concludes the proof of Theorem~$\ref{periodiclerch}$. 

\begin{remark} This, result shows that considering the periodic case does not produce new numbers. 
The special case of purely periodic generalized polylogarithmic functions corresponds to $b(z)=\prod_{i=1}^l(z^q-\alpha_i^q)$, that is $m=lq$ and roots $\zeta_l^k\alpha_i$. 
{{To prove Corollary $\ref{coro 1}$, we apply Theorem $\ref{periodiclerch}$ for $b(z):=\prod_{i=1}^m(z^{q_i}-\alpha^{q_i}_i)$ and $$w_{l}(z):=w_{i,l_i}(z)\prod_{j\neq i}(z^{q_j}-\alpha^{q_j}_j)$$ for $1\le i \le m$ and $0\le l_i \le q_i-1$.}}

It is worthwhile noting that restricting to the case where $b(z)$ has only simple roots is necessary. Indeed, if $b(z)$ does have multiple roots, an argument similar to Lemma~$\ref{elementsimple}$ shows that one gets {\it linearly dependent numbers}.
\end{remark}}
\section{Examples}
\label{ex}

{{\begin{example} \label{ex2}
Let $p$ be a prime number, $q$ a {{positive integer}}. We take $d=p$, $r_1=\cdots=r_p=10$, $m=10$ and 
$$\boldsymbol{\alpha}:=\left(1,\ldots,\dfrac{1}{10}\right), \  \boldsymbol{x}:=\left(0,\dfrac{1}{p},\ldots,\dfrac{p-1}{p}\right)\enspace.$$
Let $g,M\in \N$ with $g,M\ge 2$. Define $f_{M,g}(X)\in \Q[X]$ by
\begin{align*}
&f_{M,g}(X)=\left(2+\dfrac{1}{M}\right)X^g-\dfrac{2}{M}X^{g-1}-2X+\dfrac{2}{M}\enspace.
\end{align*}
Using the Eisenstein's irreducibility criterion of polynomial, we obtain that $f_{M,g}(X)$ is an irreducible polynomial in $\Q[X]$. 
Let $1/\beta_M$ be a root of $f_{M,g}(X)$ and put $K=\Q(\beta_M)$.
Let $w_{i,0}(z),\ldots,w_{i,q-1}(z)\in \Q[z]$ which are linearly independent over $\Q$ for $1\le i \le m$.
Let $v_0$ be an archimedean place of $K$ such that the absolute value of ${\rm{log}}\, |1/\beta_M|_{v_0}$ takes the minimal value among $({\rm{log}}\,|1/\beta_M|_v)_{v\in {{\mathfrak{M}}}^{\infty}_{\Q(\beta_M)}}$.
We apply Corollary $\ref{coro 1}$ with $q_i=q$, $\boldsymbol{\alpha},\boldsymbol{x},w_{i,l_i}(z)$ and $K$.

\bigskip

Let $M_0$ be the minimal positive integer with $|\beta^{(h)}_M|\le 2$ for $2\le h \le g$ and $M$ a positive integer with $M\ge M_0$.
Since ${\rm{den}}(\beta_M)\le 2$ and ${\rm{log}}\,|\beta_M|_v\le [K_v:\R]{\rm{log}}(2)/g$ for $v|\infty$, {{we have
the following tables showing  values of  $\log(M_0)$ such that the right-hand of the inequality below becomes positive whenever $M\ge M_0$  for each tuple of $g, p, q$:}}
\begin{align*}
V(\boldsymbol{\alpha},\boldsymbol{x},\beta_M)&\ge \dfrac{1}{{{g}}}{\rm{log}}\,|\beta_M|-\dfrac{100 pq}{{{g}}}\left[(g-1)\log(2)+\log(2520)\right]-10p\log(p)\\
&-\left[100pq\log(2)+100\left(\log(100pq+1)+100pq\log\left(\dfrac{100pq+1}{100pq}\right)\right)\right]-1000pq\enspace.
\end{align*}

\

For example, the case ${{g}}=2$ is described as follows:

$$
\begin{array}{ccccc}
  p\backslash q & 1 & 2 & 3 & 4 \\
 2 & 3158 & 5816 & 8449 & 11072 \\
 3 & 4509 & 8466 & 12398 & 16320 \\
 5 & 7192 & 13748 & 20278 & 26798 \\
 7 & 9868 & 19021 & 28150 & 37268 \\
\end{array}
$$

the case $g=3$:

$$
\begin{array}{ccccc}
 p\backslash q & 1 & 2 & 3 & 4 \\
 2 & 4427 & 8104 & 11744 & 15368 \\
 3 & 6298 & 11769 & 17202 & 22620 \\
 5 & 10013 & 19071 & 28092 & 37097 \\
 7 & 13717 & 26362 & 38969 & 51562 \\
\end{array}
$$

the case $g=4$:

$$
\begin{array}{ccccc}
 p\backslash q & 1 & 2 & 3 & 4 \\
 2 & 5695 & 10391 & 15038 & 19664 \\
 3 & 8087 & 15071 & 22006 & 28920 \\
 5 & 12834 & 24394 & 35905 & 47395 \\
 7 & 17565 & 33702 & 49789 & 65855 \\
\end{array}
$$

the case $g=5$:

$$
\begin{array}{ccccc}
 p\backslash q & 1 & 2 & 3 & 4 \\
 2 & 6964 & 12679 & 18332 & 23960 \\
 3 & 9876 & 18374 & 26809 & 35219 \\
 5 & 15655 & 29718 & 43719 & 57694 \\
 7 & 21414 & 41042 & 60608 & 80148 \\
\end{array}
$$

and the case $g=6$:

$$
\begin{array}{ccccc}
 p\backslash q & 1 & 2 & 3 & 4 \\
 2 & 8233 & 14967 & 21627 & 28256 \\
 3 & 11665 & 21676 & 31613 & 41519 \\
 5 & 18476 & 35041 & 51532 & 67992 \\
 7 & 25263 & 48382 & 71427 & 94441 \\
\end{array}
$$
demonstrated.

\end{example}}}

\noindent
{\bf Acknowledgements}

\smallskip

{{This work is
supported by JSPS KAKENHI Grant no.~18K03225, Grant  no.~21K03171 to make occasions to let the authors discuss in person, and
partly  supported by the Research Institute for Mathematical Sciences, an international joint usage
and research center located in Kyoto University.
We also deeply thank the anonymous referee for the work giving us very  precise comments as well as expert advice.
We appreciate the support by Professor Masaru Ito for calculating examples above using Mathematica.}} 

\newpage

\bibliography{}

\begin{scriptsize}
\begin{minipage}[t]{0.38\textwidth}

\noindent
Sinnou David,
\\sinnou.david@imj-prg.fr
\\Institut de Math\'ematiques \\de Jussieu-Paris Rive Gauche,
\\CNRS UMR 7586,
\\Sorbonne Universit\'{e}, 4, place Jussieu,\\75005 Paris, France\\
\& CNRS UMI 2000 Relax, \\Chennai Mathematical Institute,\\ 
 H1, SIPCOT IT Park, Siruseri, \\Kelambakkam 603103, India \\\\
\end{minipage}
\begin{minipage}[t]{0.36\textwidth}
Noriko Hirata-Kohno, 
\\hiratakohno.noriko@nihon-u.ac.jp,
(hirata@math.cst.nihon-u.ac.jp),
\\Department of Mathematics, \\College of Science \& Technology, \\Nihon University,
\\Kanda, Chiyoda, Tokyo, \\101-8308, Japan\\\\
\end{minipage}
\begin{minipage}[t]{0.35\textwidth}
Makoto Kawashima,
\\kawashima.makoto@nihon-u.ac.jp
\\Department of Liberal Arts \\and Basic Sciences, \\College of Industrial Engineering, \\Nihon University,
Izumi-chou, \\Narashino, Chiba, 275-8575, Japan
\end{minipage}

\end{scriptsize}

\end{document}